\documentclass[leqno,a4paper,12pt]{article}

\usepackage{a4wide}
\usepackage{amsmath}
\usepackage{amsfonts}
\usepackage{nicefrac}
\usepackage{mathrsfs}
\usepackage{dsfont}
\usepackage{mparhack}
\reversemarginpar

\newtheorem{theorem}{Theorem}[section]

\newtheorem{lemma}[theorem]{Lemma}
\newtheorem{proposition}[theorem]{Proposition}

\newtheorem{example}[theorem]{Example}

\numberwithin{equation}{section}

\newenvironment{proof}{{\bf
Proof:\,}}{\hspace*{\fill}\rule{1.2ex}{1.2ex}\\ }
\newenvironment{proofo}[1]{{\bf
#1:\,}}{\hspace*{\fill}\rule{1.2ex}{1.2ex}\\ }

\newcommand{\norm}[1]{\left\Vert#1\right\Vert}
\newcommand{\onorm}[1]{\Vert#1\Vert}
\newcommand{\bigset}[1]{\big\{#1\big\}}
\newcommand{\leb}{\boldsymbol{\lambda}}
\newcommand{\bignorm}[1]{\big\Vert#1\big\Vert}
\newcommand{\oabs}[1]{\vert#1\vert}
\newcommand{\scaleZ}[1]{#1^{-1}\Z^d}

\newcommand{\abs}[1]{\left\vert#1\right\vert}
\newcommand{\Bigabs}[1]{\Big\vert#1\Big\vert}
\newcommand{\bigabs}[1]{\big\vert#1\big\vert}

\newcommand{\set}[1]{\left\{#1\right\}}

\newcommand{\R}{\mathds{R}}
\newcommand{\E}{\mathds{E}}
\newcommand{\Z}{\mathds{Z}}
\newcommand{\N}{\mathds{N}}
\newcommand{\Prob}{\mathds{P}}

\newcommand{\sE}{\mathscr{E}}
\newcommand{\sX}{\mathscr{X}}

\newcommand{\cO}{\mathcal{O}}
\newcommand{\cP}{\mathcal{P}}
\newcommand{\cQ}{\mathcal{Q}}
\newcommand{\cE}{\mathcal{E}}
\newcommand{\cN}{\mathcal{N}}

\newcommand{\eps}{\varepsilon}
\newcommand{\supl}{\sup\limits}
\newcommand{\suml}{\sum\limits}
\newcommand{\liml}{\lim\limits}
\newcommand{\il}{\int\limits}
\newcommand{\iil}{\iint\limits}
\newcommand{\diag}{\mathrm{diag}}

\DeclareMathOperator{\dist}{dist}

\DeclareMathOperator{\supp}{supp}

\begin{document}

\author{Ryad Husseini\footnote{Research financed by DFG (German Science Foundation) through SFB 611} \and Moritz Kassmann\footnote{Research partially supported by DFG (German Science Foundation) through SFB 611} \footnote{corresponding author, email: {\it kassmann@iam.uni-bonn.de}}}
\title{Markov chain approximations for symmetric jump processes}

\maketitle
%

\parindent0ex
\parskip1.6ex


\begin{abstract}
  Markov chain approximations of symmetric jump processes are investigated. Tightness results and a central limit theorem are established. Moreover, given the generator of a symmetric jump process with state space $\R^d$ the approximating Markov chains are constructed explicitly. As a byproduct we obtain a definition of the Sobolev space $H^{\alpha/2}(\R^d)$, $\alpha \in (0,2)$, that is equivalent to the standard one.
\end{abstract}



\newpage

\section{Introduction}\label{Introduction}

Let $Y=(Y_t)_t$ be a time-continuous Markov chain on $\Z^d$. It is a natural question whether the sequence $(Y^n)$ of Markov chains defined by $Y^n_t = n^{-1} Y_{n^\alpha t}$, $\alpha \in (0,2]$, tends to some reasonable process for $n \to \infty$. The case $\alpha =2$ is known as diffusive scaling and leads to a diffusion process under certain assumptions on $(Y_t)_t$, see the classical Donsker's Invariance Principle of \cite{Donsker51} for the Brownian motion and chapter 11 of \cite{StrVar79} for diffusion processes in non-divergence form. In the case of symmetric processes Stroock and Zheng derive in \cite{StrZheng97} a central limit theorem for continuous-time symmetric Markov chains of bounded range. In a recent paper, Bass and Kumagai \cite{BaKu06} remove the restriction of bounded range by replacing it by a second moment condition. In both publications, the generator of the limit object is of the form $Lu(x)=\suml_{i,j=1}^d \partial_{x_i} \big( a_{i,j}(\cdot) \partial_{x_j} \big)$ and a formula is provided how the diffusion coefficient functions $a_{i,j}(\cdot)$ can be computed from the conductivities of the chain $(Y_t)$. The other direction, i.e. constructing a sequence of approximating Markov chains for a given diffusion matrix is not less important and one of the main results of \cite{StrZheng97}.

The aim of this work is to prove results analogous to ones of \cite{StrZheng97}, \cite{BaKu06} in the case where the limit object is a symmetric jump process with corresponding Dirichlet form $\big(\mathcal{E}, D(\mathcal{E})\big)$ given by
\begin{align}
\begin{split}
\mathcal{E}(f,g) &= \frac12 \iil_{(\R^d\times \R^d) \setminus \operatorname{diag}} \big(f(y)-f(x)\big)\big(g(y)-g(x)\big) k(x,y) \, dx\, dy , \quad k(x,y)=k(y,x) \,, \\
D(\mathcal{E})&= \overline{C^1_c(\R^d)}^{\cE_1}, \quad \text{ where }\cE_1(f,f)= \cE(f,f) + \|f\|_{L^2} \,,
\end{split}
\label{eq:defE}
\end{align}
and generator $L$ given by
\begin{align} Lu(x) = \liml_{\eps \to 0} \il_{|y-x|\geq \eps} \big(u(y)-u(x) \big) k(x,y) \,dy \,. \label{eq:defL}
\end{align}
Therefore, we study Markov chain approximations for a certain class
of symmetric jump processes. In  \cite{StrZheng97}, \cite{BaKu06}
the generator of the limit object is a uniformly elliptic operator.
In our situation the equivalent concept of uniform ellipticity would
be given by $k(x,y) \geq c |x-y|^{-d-\alpha}$ $\forall \, |x-y| \leq
r_0$ for some $c>0$, $r_0>0$, $\alpha \in (0,2)$. One feature of our
approach is that our central limit theorem allows for cases where
such an estimate does not hold, i.e. the limit process may be a pure
jump process which is anisotropic in some sense. The level of
anisotropy is limited since our approach uses a-priori bounds for
the modulus of continuity of the heat kernel. As discussed in
\cite{BBCK06} these bounds fail for very irregular jump measures.

There are several other contributions to the question how to
approximate Hunt processes given by Dirichlet forms, see
\cite{MRZ98}, \cite{MRS00} and the references therein. However, our
results are not covered by these works. We close the introduction by
commenting on the differences between this work and \cite{BaKu06},
\cite{StrZheng97}.
\begin{itemize}
\item The limit object in \cite{StrZheng97}  and \cite{BaKu06}  is a diffusion whereas here it is a jump process.
\item The main result of \cite{BaKu06} is a central limit theorem. In addition to such a result we establish an approximation result for a given symmetric jump process. Theorem \ref{theo:approx} should be contrasted with Theorem 3.9 from \cite{StrZheng97}.
\item Our assumption (A5) differs from (A5) of \cite{BaKu06} significantly. First, we do not assume continuity of the coefficients of the limit process. Second, we assume only $L^1_{loc}$-convergence of conductivities which is substantially less than uniform convergences on compacts.
\end{itemize}

The paper is organized as follows. In section \ref{sect:assumres} we present our assumptions and results. We provide a detailed discussion of the assumptions and some definitions and notation. Furthermore, an auxiliary result on equivalent norms on the Sobolev space $H^{\alpha/2}(\R^d)$ is proved. Sections \ref{sect:upperbounds} and \ref{sect:regularity} provide the proofs of Theorem \ref{theo:tightness-prelim} and Theorem \ref{theo:regularityheatkernels} both of which are crucial to the proof of our main results. In section \ref{sect:centrallimittheorem} we prove Theorems \ref{theo:meta-one} and \ref{theo:meta-two}. Theorem \ref{theo:approx} is proved in section \ref{sect:approx}.

{\bf Acknowledgments:} The authors express their gratitude to Professor Z.-Q. Chen and to Professor T. Kumagai for helpful comments on the presentation of the results.

\section{Assumptions and results}\label{sect:assumres}

We formulate our assumptions and results in section \ref{sect:assumform}.
In section \ref{sect:discussionofassumptions} we provide a detailed discussion of the assumptions. Section \ref{sec:sobolevnorms} is devoted to a result on equivalent norms on $H^{\alpha/2}(\R^d)$, $\alpha\in (0,2)$.  In section \ref{sect:definitions} we define and list various further objects that we deal with in this article.

We denote the counting measure by $\mu$ and the Lebesgue measure on $\R^d$ by $\leb$. In our context we also deal with the function spaces $L^2(\scaleZ{\rho}, \rho^{-d}\mu)$ where $\rho>0$. The scaling factor $\rho^{-d}$ in front of $\mu$ is natural from a geometric point of view. Write $B^\rho(x,r) := B(x,r) \cap \rho^{-1} \Z^d$ for the $r$-ball around $x$ in $\scaleZ{\rho}$. We also use the notation $\mu^\rho = \rho^{-d} \mu$. For $x \in \R^d$ we use the abbreviation $|x|_\infty = \max\limits_{i=1,\ldots,d}|x_i|$. For $x \in \R$ we write $\lceil x\rceil$ instead of $\max\{l\in \Z: l \leq x \}$.  For a point $x \in \R^d$ we denote by $[x]_n$ the element of $\scaleZ{n}$ satisfying $([x]_{n})_i=n^{-1}\lceil nx_i \rceil$ for all $i=1,\ldots,d$.

\subsection{Formulation of assumptions and results}\label{sect:assumform}
Let $(C^n)_{n \in \N}$ be a sequence of {\it conductivity functions} $C^n\colon \scaleZ{n} \times \scaleZ{n} \to [0,\infty)$. The following assumptions will be important.
\begin{itemize}
  \item[(A1)]  $C^n(x,y) = C^n(y,x)$ and $C^n(x,x)= 0$ for all $n \in \N$, $x,y \in \scaleZ{n} $.
  \item[(A2)]  There exists $\kappa_1>0$ such that for all $n \in \N$, $x\neq y$
  \begin{align*}
    C^n(x,y) \leq \kappa_1 \abs{x-y}^{-d-\alpha}.
  \end{align*}
  \item[(A3)] There exist $N_0\in\N$ and $\kappa_2>0$ with the following property: For any $n \in \N$, $x,y\in\scaleZ{n}$, $x\neq y$ there are elements $z_0^{(x,y)},\ldots, z_l^{(x,y)} \in \scaleZ{n}$, $l\leq N_0$, $z_0^{(x,y)}= x$, $z_l^{(x,y)}=y$ satisfying
    \begin{align*}
    C\big(z_i^{(x,y)}, z_{i+1}^{(x,y)}\big) \geq \kappa_2 \abs{x-y}^{-d-\alpha} \,, \quad i=0,\ldots ,l-1.
  \end{align*}
  Moreover, for any pair $(\zeta,\xi)\in \scaleZ{n} \times \scaleZ{n}$ the number of pairs $(x,y)\in\scaleZ{n}\times\scaleZ{n}$ such that $\zeta=z_k^{(x,y)}$ and $\xi=z_{k+1}^{(x,y)}$ for a $k$ is bounded by $N_0$.
\end{itemize}
For given $x,y\in\scaleZ{n}$, $x\neq y$ we call the ordered set $\{z_0^{(x,y)},\ldots, z_l^{(x,y)}\}$ above a {\it chain} and $l$ the length of the chain. The above assumptions are essential for our approach. For a mere technical reason discussed below in detail we need an additional assumption:
\begin{itemize}
  \item[(A4)] There exist $\Theta_1>0$ and $\kappa_3>0$ such that for all $x\in\scaleZ{n}$ and $r\geq \Theta_1 n^{-1}$
  \begin{align*}
    \mu \big(\bigset{y\in B^n(x, r): C^n(x,y)\geq \kappa_3 \abs{x-y}^{-d-\alpha}}\big) \geq \tfrac 56 \, \mu \big(B^n(x,r)\big)\,.
  \end{align*}
\end{itemize}
It is important for our results that the constants $\kappa_1, \kappa_2, \kappa_3, N_0, \Theta_1$ appearing in (A1) through (A4) do not depend on $n \in \N$. We associate to $C^n$ a symmetric discrete-time Markov chain $X^n=(X^n_k)_{k\in \N}$ by
\begin{align}\label{eq:defndiscretetimemarkov}
\Prob ^x (X^n_1 =y) = \frac{C^n(x,y)}{\suml_{z\in\scaleZ{n}} C^n(x,z)} \,.
\end{align}
Let $Y^n=(Y_t^n)_t$ be the symmetric time-continuous Markov chain that has the same jumps as $(X_k)$ while its holding time in the point $x$ is exponentially distributed with parameter $\sum_{z\in\scaleZ{n}} C^n(x,z)$. Note that each $Y^n$ starting in $x\in \scaleZ{n}$ corresponds to a probability measure on $D([0,\infty);\R^d)$, the space of right-continuous paths in $\R^d$ having left limits, see \cite{EtKu86}, \cite{Bil99} for properties of $D([0,\infty);\R^d)$. Our first result reads as follows:

\begin{theorem}\label{theo:meta-one}
Let $(C^n)_n$ be a sequence of conductivity functions satisfying (A1) through (A4). For $x_n \in \scaleZ{n}$, $x_n \to x \in \R^d$ the laws of $Y^n$ starting in $x_n$ are tight in $D([0, t_0];\R^d)$ for any $t_0 >0$.
\end{theorem}
For a reformulation and the proof of this result see Theorem \ref{theo:tightness}.

In order to establish a central limit theorem one needs to prescribe the behavior of $C^n$ for $n$ tending to infinity. For $x\in \scaleZ{n}$ set $\mathfrak{Q}_n(x)= \prod_{i=1}^d [x_i, x_i+1/n)$ and $\mathfrak{Q}_n=\bigcup\limits_{x\in\scaleZ{n}}\set{\mathfrak{Q}_n(x)}$.

\begin{itemize}
  \item[(A5)] There exists a measurable function $k\colon\R^d\times \R^d\to [0,\infty)$ such that for any compact subset $K$ of $\R^d\times\R^d\setminus \set{(x,x) : x\in\R^d}$ the functions $(x,y)\mapsto C^n([x]_n,[y]_n)$ converge in $L^1(K)$ to $k(\cdot,\cdot)$ for $n\to\infty$.
  \item[(B)]  There exist $M_0\in\N$ and $\Lambda_2>0$ with the following property: For any $\eps > 0$, $n \in \N$ and $\cO, \cQ \in \mathfrak{Q}_n$ there are elements $\cP_0^{(\cO,\cQ)},\ldots, \cP_l^{(\cO,\cQ)} \in \mathfrak{Q}_n$, $l\leq M_0$ with $\cP_0^{(\cO,\cQ)} = \cO$, $\cP_l^{(\cO,\cQ)}=\cQ$ and satisfying
    \begin{align*}
    \iil_{\cP_j^{(\cO,\cQ)} \times \cP_{j+1}^{(\cO,\cQ)}} k(x,y)  \mathds{1}_{\{|x-y|\geq \eps \}} \, dx dy \geq \Lambda_2 \iil_{\cO \times \cQ} \mathds{1}_{\{|x-y|\geq \eps \}} |x-y|^{-d-\alpha} \, dx dy \; \forall \, j \leq l-1 \,.
  \end{align*}
  Moreover, for any $n\in\N$ and any pair $(\mathcal{R},\mathcal{S})\in \mathfrak{Q}_n  \times \mathfrak{Q}_n$ the number of pairs $(\cO,\cQ)\in \mathfrak{Q}_n  \times \mathfrak{Q}_n$ such that $\mathcal{R}=\cP_k^{(\cO,\cQ)}$ and $\mathcal{S}=\cP_{k+1}^{(\cO,\cQ)}$ for a $k$ is bounded by $M_0$.
\end{itemize}
Again, we call the ordered set $\{\cP_0^{(\cO,\cQ)},\ldots, \cP_l^{(\cO,\cQ)}\}$ above a {\it chain} and $l$ the length of the chain. Although, to some extent, (B) is  a continuous analog of (A3) it does not follow from (A3) and (A5). Such an implication could easily be achieved by adding an additional assumption. In order not to weaken Theorem \ref{theo:meta-one} we prefer to work with (B) separately.

Here is our central limit theorem.
\begin{theorem}\label{theo:meta-two}
Let $(C^n)_n$ be a sequence of conductivity functions satisfying (A1) through (A5) and (B). Let $\sX$ be the symmetric jump process associated to the regular Dirichlet form $(\mathcal{E}, D(\mathcal{E}))$ given in (\ref{eq:defE}) and $\cN$ the properly exceptional set. Then, for $x_n \in \scaleZ{n}$, $x_n \to x \in \R^d\setminus \cN$  the laws of $Y^n$ starting in $x_n$ converge weakly in $D([0,t_0];\R^d)$ to the law of $\sX$ starting in $x$.
\end{theorem}
{\bf Remark:} The assumptions (A3) and (B) are technically involved and cover anisotropic situations. In fact, (A3) and (B) are trivially satisfied in the isotropic case, i.e. if $C^n$ satisfies (A1), (A2), (A5) and $C^n(x,y) \geq c \abs{x-y}^{-d-\alpha}$ for all $n\in \N$ and $|x-y|>n^{-1}K$ for some $K>0, c>0$. Even in  this case our theorem is still interesting and new.

It is necessary to allow for some exceptional set in Theorem \ref{theo:meta-two}. However, due to results in  \cite{ChenKumagaiHeatKernelStable} the set $\cN$ can be removed or assumed to be empty in several situations. There is no need for an exceptional set in our third result, Theorem \ref{theo:approx}. As in Theorem 3.14 of \cite{StrZheng97} we give an explicit construction of approximating Markov chains.

\begin{theorem}\label{theo:approx}
Let $k\colon\R^d\times\R^d\to (0,\infty)$ be a measurable function satisfying $k(x,y)=k(y,x)$ and
\begin{align}\label{eq:approx-assum}
  \kappa_4 \abs{x-y}^{-d-\alpha}\leq k(x,y)\leq \kappa_5 \abs{x-y}^{-d-\alpha}
\end{align}
for all $x,y\in\R^d$, $x\neq y$ with some positive constants $\kappa_4\leq\kappa_5$. Define the conductivity functions $C^n\colon \scaleZ{n} \times \scaleZ{n} \to[0,\infty)$ by
\begin{align*}
C^n(x,y) = \begin{cases} 0 \qquad &\text{ for } |x-y|_\infty \leq  n^{-1} \,, \\
n^{2d}\int\limits_{\substack{\abs{x-\xi}_\infty< n^{-1}/2\\\abs{y-\zeta}_\infty<n^{-1}/2}} k(\xi,\zeta) d\xi\,d\zeta &\text{ for } |x-y|_\infty \geq 2 n^{-1}
\end{cases}
\end{align*}
Let $\sX$ be the Hunt process corresponding to the Dirichlet form $\big(\mathcal{E}, D(\mathcal{E})\big)$ given by (\ref{eq:defE}).
Then the sequence of processes corresponding to $C^n$ converges in the sense of Theorem \ref{theo:meta-two} to $\sX$ for any starting point $x\in \R^d$.
\end{theorem}

{\bf Remark:} As becomes clear by the discussion below (A4) allows for quite general cases of sequences $C^n$. In addition to it, in light of Lemma \ref{lem:necessaryA2A3} and Lemma \ref{lem:ondiagonalupperboundlargescale} it is very likely that (A4) can be dropped. This would imply the possibility to weaken the lower bound on $k$ assumed in (\ref{eq:approx-assum}) substantially.

\subsection{Discussion of assumptions}\label{sect:discussionofassumptions}

We illustrate assumptions (A1) through (A5) introduced in Section \ref{sect:assumform}. First, let us look at (A1) through (A4). If, for a fixed scale $n\in\N$, $C^n\colon \scaleZ{n}\times\scaleZ{n}\to \R^+$ satisfies $(A1)$ through $(A4)$ then the same holds for the conductivity function $C\colon \Z^d\times\Z^d\to \R^+$ with $C(x,y)=n^{-d-\alpha} C(n^{-1} x, n^{-1} y)$ for $x,y \in\Z^d$ with the same constants $d, \kappa_1,\kappa_2, N_0$ where the chains in (A3) have to be scaled in an obvious way. Naturally, the interaction radius $\Theta_1$ in (A4) is multiplied by $n$. Since $C$ is the appropriate conductivity function corresponding to the process $Y^n$ scaled on $\Z^d$ in the obvious (''$\alpha$-stable'') way, it is sufficient to understand (A1) through (A4) for a single conductivity function $C$ on $\Z^d\times\Z^d$. In addition, the main results of Section 3 and 4 are scale-invariant, i.e. the constants appearing are scale-invariant, and depend only on the constants in the assumptions. Therefore it again suffices to prove them for a fixed conductivity function.

(A3) and (A4) are stable under perturbations of the conductivity function near the diagonal. This reflects the fact that in our central limit theorem jumps smaller than a fixed $R>0$ have no influence on the limit process. Let $C$ and $\widetilde{C}$ be conductivity functions  such that there exists $R>0$ with $\widetilde{C}(x,y)=C(x,y)$ for $\abs{x-y}\geq R$. If $C$ satisfies one of the assumptions (A3), (A4) then the same assumption also holds for $\widetilde{C}(x,y)$ with the same constants. On the other hand, if $C^n, \widetilde {C}^n\colon \scaleZ{n}\times\scaleZ{n}\to\R^+$ are to sequences of conductivities with $\widetilde{C}^n (x,y) = C^n(x,y)$ whenever $\abs{x-y}\geq Rn^{-1}$ and if $(C^n)$ satisfies (A5), then (A5) also holds for $(\widetilde{C}^n)$ with the same limit function $k$.

(A1) implies the process to be symmetric while (A2) bounds $C(x,y)$ from above by the conductivities of a rotationally symmetric  $\alpha$-stable Markov chain on $\Z^d$. (A2) gives in particular
\begin{align*}
  \sup_{x\in\Z^d} \sum_{y\in\Z^d} C(x,y)<\infty.
\end{align*}

(A3) is much more technical. It implies a certain kind of irreducibility of the associated Markov chain. Additionally, it takes into account the highly non-local nature of our objects. Roughly it says that every two points $x,y$ can be connected by chaining together at the utmost $N_0$ jumps where the probability of each jump is bounded from below by a constant multiple of $\abs{x-y}^{-d-\alpha}$ while at the same time one has enough of these connecting jumps.

(A2) and (A3) imply together $\bigabs{z_i^{(x,y)} - z_{i+1}^{(x,y)}}\leq (\kappa_1/\kappa_2)^{1/(d+\alpha)}\abs{x-y}$. This leads us to the following necessary condition for (A2) and (A3).
\begin{lemma}\label{lem:necessaryA2A3}
  Assume (A2) and (A3). Then there exist $\gamma\in (0,1)$, $\Theta_1>0$ and $\kappa_3>0$ depending only on $\kappa_1$, $\kappa_2$, $N_0$, $d$ and $\alpha$ such that for all $x\in\Z^d$, $r>\Theta_1$
  \begin{align}\label{eq:necessaryA2A3}
    \mu\big(\{y\in B^1(x,r) : C(x,y)\geq \kappa_3\abs{x-y}^{-d-\alpha}\}\big )\geq \gamma\mu\big(B^1(x,r)\big) .
  \end{align}

  In particular, if the conductivities are stationary, i.e. $C(x,y) = \widetilde{C}(x-y)$ then
  \begin{align*}
    \mu\big(\{h \in B^1(0,r): \widetilde{C}(h)\geq \kappa_3\oabs{h}^{-d-\alpha}\}\big)\geq \gamma\mu(B^1(0,r)).
  \end{align*}
\end{lemma}
\begin{proof}
  First notice that (A2) and (A3) imply the existence of $c_1= c_1(\kappa_1, \kappa_2, d, \alpha, N_0)\geq 1$ such that for all $l\leq N_0$, $\xi,\zeta\in\Z^d$
  \begin{align*}
    \bigabs{z_l^{(\xi,\zeta)}-\xi}\leq c_1\abs{\xi-\zeta}.
  \end{align*}
  Assume $r$ large enough and $x\in\Z^d$. Then $M=\bigset{z^{(x,y)}_1 \in\Z^d\colon y\in B^1(x,r/c_1)}\subset B^1(x,r)$. By the second part of (A3)
  \begin{align}
   \mu(M) \geq \lceil\mu\big(B^1(x,r/c_1)\big)/N_0\rceil\geq \frac{c_2(d)}{c^d_1N_0} \mu\big(B^1(x,r)\big) \,.
   \end{align}
Next, set
  \begin{align*}
    \widetilde{M}=M \setminus B^1\Big(x,\big(\frac{c_2}{2c^d_1N_0}\big)^{1/d} r\Big).
  \end{align*}
   Trivially, $\mu(\widetilde{M})\geq c_3 \mu\big(B^1(x,r)\big)$ where $c_3=\frac{c_2}{4c^d_1N_0}$ depends on all constants that appeared so far. Assume $r\geq \big(\frac{c_2}{2c^d_1N_0}\big)^{-1/d}$ and $z\in \widetilde{M}$. Then there is $y \in B(x,r/c_1)$ with $z=z^{(x,y)}_1\in \widetilde{M}$ and
  \begin{align*}
    C\big(x,z^{(x,y)}_1\big)&\geq \kappa_2\abs{x-y}^{-d-\alpha}\geq \kappa_2c^{d+\alpha}_1r^{-d-\alpha}\geq\kappa_2\big(\frac{c_2}{2N_0}\big)^{(d+\alpha)/d}\bigabs{x-z^{(x,y)}_1}^{-d-\alpha}.
  \end{align*}
Setting $\kappa_3 = \kappa_2\big(\frac{c_2}{2N_0}\big)^{(d+\alpha)/d}$ and $\gamma = c_3$ the set $\{y\in B^1(x,r) : C(x,y)\geq \kappa_3\abs{x-y}^{-d-\alpha}\}$ contains $\widetilde{M}$ and satisfies \eqref{eq:necessaryA2A3}.
\end{proof}

Lemma \ref{lem:necessaryA2A3} immediately implies that under our assumptions a second moment condition as in \cite{BaKu06} cannot hold. But (\ref{eq:necessaryA2A3}) is not sufficient neither for (A2) nor for (A3). Take for instance $d=1$ and $C(x,y)= \abs{x-y}^{-1-\alpha}$ if $x\neq y$, $x-y\in 2\Z$ and $C(x,y)=0$ else. Then clearly (A3) is not satisfied.

Let us now provide some examples of conductivity functions satisfying our assumptions. If $C(x,y) \abs{x-y}^{d+\alpha}$ stays bounded between two positive constants then $C$ satisfies (A2), (A3) and (A4). Hence all cases of \cite{BaLe02b} are covered by our conditions. In addition, our assumptions allow for cases where there are no jumps in the direction of certain cones.
\begin{example}\label{ex:discretekegel}
  Let $V:=\set{(x_1,x_2)\in\Z^2 : \abs{x_2}\leq \gamma\abs{x_1}}$, $\gamma>0$ be a double-cone in $\Z^d$. Set
  \begin{align}\label{example:doublecone}
    C(x,y) = \mathds{1}_V (x-y) g(x,y) \abs{x-y}^{-d-\alpha}
  \end{align}
  for $x,y\in\Z^2$, $x\neq y$ where $g\colon \Z^d \times \Z^d \to [a,b]$ for some $0<a<b$ is a measurable, symmetric function. Then these conductivities satisfy (A2) and (A3). If $\gamma$ is large enough (A4) holds, too.

  In fact, if $x=(x_1, x_2)$ and $y=(y_1,y_2)$ set $z^{xy} = (x_1+\lceil(2+\gamma^{-1})\abs{x-y}\rceil,x_2)$. Then
  \begin{align*}
    x-z^{xy}&=(\lceil(2+\gamma^{-1})\abs{x-y}\rceil, 0)\in V,\\
    \abs{x-z^{xy}}&\leq (3+\gamma^{-1})\abs{x-y},\displaybreak[0]\\
    z^{xy}-y &=(x_1-y_1+\lceil(2+\gamma^{-1})\abs{x-y}\rceil, x_2-y_2),\displaybreak[0]\\
    \abs{z^{xy}-y} &\leq \abs{x-y} +\sqrt{2}\lceil(2+\gamma^{-1})\abs{x-y}\rceil\leq 2 (4+\gamma^{-1})\abs{x-y}.\displaybreak[0]
  \end{align*}
  Finally, $z^{xy}-y\in V$ by
  \begin{align*}
    \gamma\abs{x_1-y_1 + \lceil(2+\gamma^{-1})\abs{x-y}\rceil}&\geq \gamma\lceil(1+\gamma^{-1})\abs{x-y}\rceil\geq \abs{x-y}\geq \abs{x_2-y_2}.
  \end{align*}
\end{example}

\begin{example}\label{ex:contkegel}
Define $C$ as in example \ref{ex:discretekegel} with $\gamma$ large enough, say $\gamma > 7/8$, and $g\equiv 1$. \\ Set $C^n(x,y)=n^{d+\alpha} C(nx, ny)$,$n \in \N$, $x,y\in \scaleZ{n}$. Then $(C^n)_n$ satisfies (A1) through (A5) with
\begin{align}
    k(x,y) = \mathds{1}_V (x-y) \abs{x-y}^{-d-\alpha}, \quad V:=\set{(x_1,x_2)\in\R^2 : \abs{x_2}\leq \gamma\abs{x_1}} \,.
\end{align}
\end{example}

Note also the following counterexample:
\begin{example}
  Set $V:= \set{0}\times \Z \cup \Z \times \set{0}$ and $C(x,y) = \mathds{1}_V (x-y) \abs{x-y}^{-d-\alpha}$. Then these conductivities do not satisfy (A2) and (A3) since the necessary condition of Lemma \ref{lem:necessaryA2A3} does not hold.
\end{example}

Finally, let us give the most obvious example of a conductivity function $C$ satisfying (A1) through (A5) .
\begin{example}
  Fix $\alpha\in(0,2)$. Then the conductivity functions \[C^n(x,y)  = \frac{\alpha \Gamma(\frac{d+\alpha}{2})}{2^{1-\alpha} \pi^{d/2}\Gamma(1-\alpha/2)}\abs{x-y}^{-d-\alpha},\quad x,y\in\scaleZ{n}\] satisfy (A1) through (A5). The limit process $\sX$ in the sense of Theorem \ref{theo:meta-two} is the well-known rotationally invariant $\alpha$-stable process. The properly exceptional set $\cN$ is empty.
\end{example}

\subsection{Equivalent norms on $H^{\alpha/2}(\R^d)$}\label{sec:sobolevnorms}

So far we have concentrated on a discussion of (A1) through (A5). Let us now look at (B). Let $k\colon\R^d \times \R^d\to [0,\infty)$ be a measurable function satisfying $k(x,y) = k(y,x)$ and $k(x,y) \leq \Lambda_1 \abs{x-y}^{-d-\alpha}$ for almost all $(x,y)$ with $x\ne y$ and some $\Lambda_1>0$. In light of (A1), (A2) and (A5) this is the structure of kernels appearing in the limit $n \to \infty$. Under assumptions (A1) through (A5) there still can be large oscillations of $C^n$ in the following sense. Fix two sequences of $x^i_n,y^i_n \in \scaleZ{n}$, $i\in\{1,2\}$, with $|x^1_n - x^2_n| \to 0$ and $|y^1_n - y^2_n| \to 0$. Then chains connecting $x^1_n$ and $y^1_n$ can be very far apart from chains connecting $x^2_n$ and $y^2_n$, no matter how large $n$ is. Assumption (B) guarantees that this phenomenon can be avoided by choosing appropriate chains. Having studied (A1) through (A3) it should be clear how to construct examples of kernels $k$ satisfying (B). For instance, the kernel $k$ constructed in example \ref{ex:contkegel} satisfies (B).

Let us show that (B) is a natural assumption. Assuming (B) we show that $D(\mathcal{E})$ from (\ref{eq:defE}) equals $H^{\alpha/2}(\R^d)$, i.e. (B) determines a class of equivalent norms on $H^{\alpha/2}(\R^d)$, $\alpha \in (0,2)$. One standard definition of this function space is
\[ H^{\alpha/2}(\R^d) = \{ f\in L^2(\R^d); \|f\|_{H^{\alpha/2}} := \sqrt{\|f\|_{L^2} + \cE_\alpha(f,f)} < \infty \} \,. \]
We show how it is possible to replace $|x-y|^{-d-\alpha}$ in the definition of $\cE_\alpha(f,f)$ by some anisotropic kernel $k(x,y)$ without changing the function space. Such a result does not seem to be established in the literature on function spaces.

\begin{theorem}\label{theo:equivnorms}
Let $k\colon\R^d \times \R^d\to [0,\infty)$ be a measurable function satisfying $k(x,y) = k(y,x)$ and $k(x,y) \leq \Lambda_1 \abs{x-y}^{-d-\alpha}$ for almost all $(x,y)$ with $x\ne y$ and some $\Lambda_1>0$, $\alpha \in (0,2)$. Assume that $k$ satisfies (B). Then there are two positive constants $c_0, c_1$ such that
\begin{align}\label{eq:formscomp}
&c_0 \cE_\alpha(f,f)  \leq \cE(f,f)  \leq c_1 \cE_\alpha(f,f) \quad \forall \, f\in C^1_c(\R^d) \,.
\end{align}
Hence, the regular Dirichlet form $\big(\mathcal{E}, D(\mathcal{E})\big)$ of (\ref{eq:defE}) satisfies under (B) \[ D(\mathcal{E})=\overline{C^1_c(\R^d)}^{\cE_1} = H^{\alpha/2}(\R^d) \,. \]
\end{theorem}

\begin{proof}
The second estimate in (\ref{eq:formscomp}) follows trivially from the upper bound of $k$. In order to establish the first one note
\[ \cE_\alpha(f,f) = \lim\limits_{\eps\to 0} \iil_{\R^d\times \R^d} \big(f(y)-f(x)\big)^2 \mathds{1}_{\{|x-y|\geq \eps \}} |x-y|^{-d-\alpha} \, dx\, dy \]
Denote by $z_{\cQ}$ the center point of a given cube $\cQ\in\mathfrak{Q}_n$. For simplicity we assume that the length of each chain is equal to $M_0$.
Assumption (B) gets involved in the following way:
{\allowdisplaybreaks
\begin{align*}
&\iil_{\R^d\times \R^d} \big(f(y)-f(x)\big)^2 \mathds{1}_{\{|x-y|\geq \eps \}} |x-y|^{-d-\alpha} \, dx \,dy \\
&= \lim\limits_{n\to\infty} \sum\limits_{\cO \in \mathfrak{Q}_n} \sum\limits_{\cQ \in \mathfrak{Q}_n} \big(f(z_{\cO})-f(z_{\cQ})\big)^2 \iil_{\cO\times \cQ} \mathds{1}_{\{|x-y|\geq \eps \}} |x-y|^{-d-\alpha} \, dx \,dy \\
&\leq M_0 \lim\limits_{n\to\infty} \sum\limits_{\cO \in \mathfrak{Q}_n} \sum\limits_{\cQ \in \mathfrak{Q}_n} \suml_{j=0}^{M_0-1} \big(f(z_{\cP_{j+1}^{(\cO,\cQ)}})-f(z_{\cP_j^{(\cO,\cQ)}})\big)^2 \iil_{\cO\times \cQ} \mathds{1}_{\{|x-y|\geq \eps \}} |x-y|^{-d-\alpha} \, dx \,dy \\
&\leq \frac{M_0}{\Lambda_2} \lim\limits_{n\to\infty} \sum\limits_{\cO \in \mathfrak{Q}_n} \sum\limits_{\cQ \in \mathfrak{Q}_n} \suml_{j=0}^{M_0-1} \big(f(z_{\cP_{j+1}^{(\cO,\cQ)}})-f(z_{\cP_j^{(\cO,\cQ)}})\big)^2 \iil_{\cP_{j+1}^{(\cO,\cQ)} \times \cP_{j}^{(\cO,\cQ)}} \mathds{1}_{\{|x-y|\geq \eps \}} k(x,y) \, dx \,dy \\
&\leq \frac{(M_0)^2}{\Lambda_2} \lim\limits_{n\to\infty} \sum\limits_{\cO \in \mathfrak{Q}_n} \sum\limits_{\cQ \in \mathfrak{Q}_n} \max\limits_{j=0,\ldots,M_0-1} \big(f(z_{\cP_{j+1}^{(\cO,\cQ)}})-f(z_{\cP_j^{(\cO,\cQ)}})\big)^2 \\
&\qquad \times \iil_{\cP_{j+1}^{(\cO,\cQ)} \times \cP_{j}^{(\cO,\cQ)}} \mathds{1}_{\{|x-y|\geq \eps \}} k(x,y) \, dx\, dy \\
&\leq \frac{(M_0)^3}{\Lambda_2} \lim\limits_{n\to\infty} \sum\limits_{\cO \in \mathfrak{Q}_n} \sum\limits_{\cQ \in \mathfrak{Q}_n} \big(f(z_{\cO})-f(z_{\cQ})\big)^2 \iil_{\cO \times \cQ} \mathds{1}_{\{|x-y|\geq \eps \}} k(x,y) \, dx \,dy \\
&= \frac{(M_0)^3}{\Lambda_2} \iil_{(\R^d\times \R^d) \setminus \operatorname{diag}} \big(f(y)-f(x)\big)^2 \mathds{1}_{\{|x-y|\geq \eps \}} k(x,y) \, dx \,dy \,.
\end{align*}
}
\end{proof}

\subsection{Further definitions and notation}\label{sect:definitions}

If $\sX$ is a stochastic process and $\Omega$ a Borel set write $\tau(\Omega; \sX)$ resp. $\sigma(\Omega; \sX)$ for the first time the process exits resp. enters $\Omega$ where we omit $\sX$ if there is no danger of confusion.

Let $X=(X_k)_k$ be the symmetric discrete-time Markov chain associated to $C$ by \eqref{eq:defndiscretetimemarkov}. A symmetric time-continuous Markov chain $Y=(Y_t)_t$ having the same jumps as $X$ can be constructed as follows: Take a family $(T_{x, j})_{x\in \Z^d, j\in\N}$ of independent random variables, independent also of $X$, such that $T_{x,j}$ is exponentially distributed with parameter $C_x$ and set $T_{x,0}\equiv 0$.
Set $Y_t=X_n$ for $t \in [\sum_{j=0}^n T_{X_j, j}, \sum_{j=0}^{n+1} T_{X_j, j})$. Note that in \cite{BaLe02b} the holding times of the time-continuous process are exponentially distributed with parameter $1$ leading to different generators and Dirichlet forms. For more details on Markov chains we refer the reader to \cite{NorrisMarkov}. (A2) and (A3) give uniform bounds on the expected holding times of $Y$. The process $Y$ corresponds to the Dirichlet form
\begin{align*}
  \mathcal{E} (f, f) &= \frac 12  \sum_{x,y\in\Z^d} \big(f(x)- f(y)\big)^2 C(x,y)\,.
\end{align*}
with domain $D(\mathcal{E})=L^2(\Z^d, \mu)$. We derive properties of $Y$ in the next section by comparing $\big(\mathcal{E}, L^2(\Z^d, \mu)\big)$ to the Dirichlet form $\big(\mathcal{E}_\alpha, L^2(\Z^d, \mu)\big)$ of a rotational invariant $\alpha$-stable process $V$ in $\Z^d$, i.e. $\mathcal{E}_\alpha$ is defined by
\begin{align*}
  \mathcal{E}_\alpha (f,f) = \frac 12 \sum_{x,y\in\Z^d}\big(f(x)-f(y)\big)^2 \abs{x-y}^{-d-\alpha}, \quad \alpha\in (0,2) \,.
\end{align*}
By (A2) $\mathcal{E} (f,f) \leq \kappa_1 \mathcal{E}_\alpha(f,f)$ for all $f\in L^2(\Z^d, \mu)$.


Define a family $(Y^\rho)_{\rho>0}$ of symmetric time-continuous Markov chains on $\scaleZ{\rho}$ by $Y_t^\rho = \rho^{-1}Y_{\rho^\alpha t}$. $Y^\rho$ corresponds to the Dirichlet form $\big(\mathcal{E}^\rho,L^2(\scaleZ{\rho}, \rho^{-d} \mu)\big)$ defined by
\begin{align}\label{eq:fullscaleddirichletform}
  \mathcal{E}^\rho (f,f) = \frac 12 \sum_{x,y\in\scaleZ{\rho}} \big(f(x) - f(y)\big)^2 C^{\rho}(x,y)  \rho^{-2d} \quad \text{ with } C^{\rho}(x,y) = \rho^{d+\alpha}C(\rho x,\rho y) \,.
\end{align}
{\bf Remark:} Note that we abuse our own notation here. The above definition of the family $(C^\rho)_{\rho>0}$ does not correspond correctly to our sequence of conductivity functions $(C^n)_{n\in\N}$ defined in the introduction. To be precise: Given an arbitrary sequence  $(C^n)_{n\in\N}$ in the sense of the introduction there might be no conductivity function $C\colon \Z^d\times \Z^d \to [0,\infty)$ such that $C^\rho=C^n$ for $\forall \, \rho=n$. Nevertheless, we use $C^\rho$ in the sense above and $C^n$ in the sense of assumptions (A1) through (A5). This remark carries over to the definition on the family $(Y^n)_{n \in \N}$. We use the symbol $Y^n$ for the time continuous process corresponding to the conductivity function $C^n$. That is, the family $(Y^\rho)_{\rho > 0}$ is determined by a single conductivity function $C\colon \Z^d\times \Z^d \to [0,\infty)$ whereas the family $(Y^n)_{n \in \N}$ depends on the whole sequence $(C^n)_{n\in \N}$. $Y^n$ corresponds to the Dirichlet form $\big(\mathcal{E}^n,L^2(\scaleZ{n}, n^{-d} \mu)\big)$ defined by
\begin{align}\label{eq:defEn}
  \mathcal{E}^n (f,f) = \frac 12 \sum_{x,y\in\scaleZ{n}} \big(f(x) - f(y)\big)^2 C^{n}(x,y) n^{-2d}  \,.
\end{align}

Note that (A1) through (A4) are stable in the following sense. If one fixed conductivity function $C$ satisfies (A1) ((A2), (A3) resp.) then the assumption holds true for the family $(C^\rho)_\rho$ with the constants independent of $\rho$. On the other hand, if (A4) is true for $C$ the conclusion of (A4) holds for any $C^\rho$ with the same $\gamma$ and $\kappa_3$ whenever $r\geq \Theta_1\rho^{-1}$.

Scaling as above implies a relation between the heat kernel $p_{Y^\rho}$ of $Y^\rho$ with respect to $\rho^{-d}\mu$ and the heat kernel $p_Y$ of $Y$. Note that, by regarding the heat kernel of the scaled process with respect to $\rho^{-d}\mu$, $p_{V^\rho} (t,x,y)$ is not anymore the probability that the process starting in $x$ is at time $t$ in $y$ but $\rho^{-d}$ times this probability. One has \begin{align}\label{eq:salingheatkernel}
  p_{Y^\rho} (t,x,y) = \rho^d p_Y(\rho^\alpha t, \rho x, \rho y).
\end{align}

Let $Y^{\rho,\lambda}$ be the process $Y^\rho$ with all jumps bigger than $\lambda$ removed. $Y^{\rho,\lambda}$ corresponds to the Dirichlet form $\big(\mathcal{E}^{\rho,\lambda},L^2(\scaleZ{\rho}, \rho^{-d} \mu)\big)$ defined by
\begin{align*}
  \mathcal{E}^{\rho, \lambda} (f,f) &=\frac 12 \sum_{\substack{x,y\in\scaleZ{\rho}\\ \abs{x-y}\leq \lambda}} \big(f(x) - f(y)\big)^2 C^{\rho}(x,y)  \rho^{-2d}.
\end{align*}

Finally, let $V^{\rho,1}$ be the process associated to
\begin{align*}
  \mathcal{E}_\alpha^{\rho,1} (f,f) = \frac 12 \sum_{\substack{x,y\in \scaleZ{\rho}\\\abs{x-y}\leq 1}} \big(f(x) - f(y)\big)^2 \abs{x-y}^{-d-\alpha} \rho^{-2d}.
\end{align*}

Let us finish this section with an overview over all processes which we have introduced so far:

\begin{minipage}[t]{0.1\textwidth}
\hskip2cm
\end{minipage}
\begin{minipage}[t]{0.8\textwidth}
\begin{enumerate}
  \item[$X=(X_k)$:]  time-discrete Markov chain on $\Z^d$ corresponding to the conductivity function $C$.
  \item[$Y=(Y_t)$:]  time-continuous Markov chain on $\Z^d$ with the same jumps as $(X_n)$; its Dirichlet form is $\mathcal{E}$.
  \item[$Y^\rho=(Y_t^\rho)$:]  scaled version of the process $(Y_t)$; it corresponds to $C^\rho$; its state space is $\scaleZ{\rho}$; its Dirichlet form is $\mathcal{E}^\rho$.
  \item[$Y^n=(Y_t^n)$:] time-continuous process corresponding to $C^n$; its state space is $\scaleZ{n}$; its Dirichlet form is $\mathcal{E}^n$ defined in (\ref{eq:defEn}).
  \item[$\sX=(\sX_t)$:] limit of $Y^n$ for $n \to \infty$; its state space is $\R^d$; corresponds to Dirichlet form $\big(\cE,D(\cE)\big)$ defined in (\ref{eq:defE}).
  \item[$Y^{\rho, \lambda}=(Y_t^{\rho, \lambda})$:] equals the process $(Y_t^\rho)$ but with jumps greater than $\lambda$ removed; its Dirichlet form is $\mathcal{E}^{\rho,\lambda}$.
  \item[$V=(V_t)$:] rotational-invariant $\alpha$-stable process on $\Z^d$; its Dirichlet form is $\mathcal{E}_\alpha$.
  \item[$V^{\rho,\lambda}=(V^{\rho,\lambda}_t)$:] scaled version of the process $V$  with jumps bigger than $\lambda$ removed; its state space is $\scaleZ{\rho}$; its Dirichlet form is $\mathcal{E}_\alpha^{\rho,\lambda}$.
\end{enumerate}
\end{minipage}

\section{Upper bounds for exit times and the heat kernel}\label{sect:upperbounds}

The aim of this section is to establish upper bounds on the heat kernel of the processes $Y^{\rho, \lambda}$ independent of $\rho\geq 1$. This is done in section \ref{sect:upperboundshk}. The results are then applied in section \ref{sect:tightness} in order to establish Theorem \ref{theo:tightness-prelim} which is the key ingredient needed to show tightness of the family $(Y^n)_{n \in \N}$. The techniques used in this section are borrowed from \cite{CaKuStr87}, \cite{BaLe02b}, and \cite{ChenKumagaiHeatKernelStable}.

\subsection{Upper bounds on the heat kernel}\label{sect:upperboundshk}

The main tools in this section are techniques worked out in \cite{CaKuStr87}. Lemma \ref{lem:comparabilitydirichletforms} is new and, together with Lemma \ref{lem:necessaryA2A3}, we consider it  important for the further development of anisotropic jump processes and Markov chains.
\begin{lemma}\label{lem:comparabilitydirichletforms}
  Assume (A1), (A2) and (A3). Then there exist $c>0$ depending only on $N_0$, $\kappa_2$ and $\Theta_2\geq 1$ depending on $\kappa_1, \kappa_2, N_0, d, \alpha$ such that for all $f\in L^2(\Z^d, \mu)$, $\rho > 0$, $\lambda > 0$
  \begin{align*}
    \mathcal{E^{\rho,\lambda}_\alpha}(f,f)\leq c\mathcal{E}^{\rho,\lambda{\Theta_2}}  (f,f) \quad \text{ and in particular } \mathcal{E^{\rho}_\alpha}(f,f)\leq c\mathcal{E}^{\rho}  (f,f) \,.
  \end{align*}
\end{lemma}
\begin{proof}
  Let $(z^{(x,y)}_l)$ be the chains associated to $C$ by (A3) now scaled on $\scaleZ{\rho}$. Note that (A2) and (A3) together imply $|z_{l-1}^{(x,y)} - z_{l}^{(x,y)}| \leq \Theta_2 |x-y|$ with $\Theta_2=\Theta_2(\kappa_1, \kappa_2, N_0, d, \alpha)$ for any chain in the sense of (A3), any pair $(x,y)$ and any $l$. For notational convenience we assume the length of all chains to be equal to $N_0$. Then
\begin{align*}
    \mathcal{E}_\alpha^{\rho,\lambda}(f,f)&=\frac12 \sum_{\substack{x,y \in\scaleZ{\rho}\\\abs{x-y}\leq \lambda}} \big(f(x) -f(y)\big) ^2 \abs{x-y}^{-d-\alpha} \rho^{-2d}\displaybreak[0]\\
    &\leq N_0 \sum_{\substack{x,y\in\scaleZ{\rho}\\\abs{x-y}\leq \lambda }} \sum_{l=1}^{N_0} \big(f(z_{l-1}^{(x,y)}) - f(z_l^{(x,y)})\big)^2 \abs{x-y}^{-d-\alpha}\displaybreak[0]\rho^{-2d}\\
    &\leq N_0 (\kappa_2)^{-1}\sum_{\substack{x,y\in\scaleZ{\rho}\\\abs{x-y}\leq \lambda}}  \sum_{l=1}^{N_0} \big(f(z_{l-1}^{(x,y)}) - f(z_l^{(x,y)})\big)^2 C^\rho(z_{l-1}^{(x,y)}, z_{l}^{(x,y)})\rho^{-2d}\displaybreak[0]\\
    &\leq N_0^2 (\kappa_2)^{-1}\sum_{\substack{x,y\in\scaleZ{\rho}\\\abs{x-y}\leq \lambda}}  \rho^{-2d} \max_{l=1,\ldots,N_0} \Big\{ \big(f(z_{l-1}^{(x,y)}) - f(z_l^{(x,y)})\big)^2  C^\rho(z_{l-1}^{(x,y)}, z_{l}^{(x,y)})\Big\} \displaybreak[0]\\
    &\leq (N_0)^3(\kappa_2)^{-1}\mathcal{E}^{\rho,\lambda{\Theta_2}}(f,f).\displaybreak[0]
  \end{align*}
For the last inequality we use the fact that every term of the sum on the left appears at least once in the sum on the right hand side. By the second part of (A3) this happens at most $N_0$ times.
\end{proof}

Recall the following on-diagonal estimate for the truncated $\alpha$-stable process on $\scaleZ{\rho}$ given in \cite{BaLe02b}:
\begin{proposition}\label{prop:upperboundtruncatedlevy}
There exists $c>0$ independent of $\rho\geq 1$ such that
\begin{align*}
  p_{V^{\rho,1}} (t,x,y) &\leq \begin{cases}
    ct^{-d/\alpha}  &\text{ for } t\in (0,1] \,,\\
    ct^{-d/2}       &\text{ for } t>1 \, ,
  \end{cases}
\end{align*}
  i.e. in particular $p_{V^{\rho,1}} (t,x,y) \leq c t^{-d/\alpha} e^t$  for all $t > 0$.
\end{proposition}
These estimates lead almost directly to upper bounds for the heat kernel of $Y^{\rho,\lambda}$.
\begin{lemma}\label{lem:ondiagonalupperboundlargescale}
  Assume (A1), (A2) and (A3). Let $\Theta_2$ be the constant arising in the statement of Lemma \ref{lem:comparabilitydirichletforms}. Then there exists $c>0$  such that for all $\rho\geq1$ and $\lambda\geq \Theta_2$
  \begin{align*}
    p_{Y^{\rho, \lambda}} (t,x,y) \leq \begin{cases}
      ct^{-d/\alpha}  &\text{for $t\leq 1$},\\
      ct^{-d/2}       &\text{for $t>1$.}
    \end{cases}
  \end{align*}
\end{lemma}
\begin{proof}
Applying Theorem 2.1 and Theorem 2.9 of \cite{CaKuStr87} to the assertions of Proposition \ref{prop:upperboundtruncatedlevy} one obtains
\begin{align*}
    \norm{f}^{2+ 2\alpha/d}_{L^2}  & \leq c_1\big(\norm {f}_{L^2}^2+\mathcal{E}_\alpha^{\rho,1}(f,f)\big)\norm{f}^{2\alpha/d}_{L^1} \,, \\
    \norm{f}^{2+ 4/d}_{L^2}  &\leq c_2\big(\norm {f}_{L^2}^2+\mathcal{E}_\alpha^{\rho,1}(f,f)\big)\norm{f}^{4/d}_{L^1} \,.
  \end{align*}
Estimating $\mathcal{E}_\alpha^{\rho,1}(f,f)$ from above with the help of Lemma \ref{lem:comparabilitydirichletforms} and applying afterwards the converse parts of Theorem 2.1 and Theorem 2.9 of \cite{CaKuStr87} the desired result follows.
%
\end{proof}

In a similar fashion, one can use the scaling property \eqref{eq:salingheatkernel}, Lemma \ref{lem:comparabilitydirichletforms} and the upper bounds on the rotational alpha-stable process on $\Z^d$ (see Proposition 4.2 in \cite{BaLe02b}) to obtain directly upper bounds on the heat kernel of $(Y^\rho_t)$:
\begin{lemma}\label{lem:ondiagonalupperboundallscales}
  Assume (A2) and (A3). Then there exists a constant $c>0$ independent of $\rho>0$ such that
  \begin{align*}
    p_{Y^\rho}(t,x,y) \leq c t^{-d/\alpha} \quad \forall \, t>0,  \forall \, x,y \in \scaleZ{\rho} \,.
  \end{align*}
  In fact, $c$ only depends on $d$, $\alpha$ and the constants $N_0$ and $\kappa_2$ appearing in (A3).
\end{lemma}

Off-diagonal estimates on the heat kernels of the truncated process can be obtained by Davies' method, see for example {\S}3 in \cite{CaKuStr87}.
\begin{lemma}\label{lem:offdiagonalboundtruncated}
  Assume (A1), (A2) and (A3). For any $\lambda\geq \Theta_2$ there exists $c>0$ such that for all $\rho\geq 1$, $x,y\in \scaleZ{\rho}$ and $t\in (0,1]$
  \begin{align*}
    p_{Y^{\rho,\lambda}} (t,x,y) &\leq c t^{-d/\alpha}e^{-\abs{x-y}}.
  \end{align*}
\end{lemma}
\begin{proof}
  Fix $x,y\in\scaleZ{\rho}$, $x\neq y$. Define
  \begin{align*}
    \Gamma (f,f) (x) &= \sum_{\substack{y\in\scaleZ{\rho}\\\abs{x-y}\leq \lambda}} \big(f(x)-f(y)\big)^2 C^\rho(x,y) \rho^{-d},\\
    \Lambda^2(\psi)&=\max\set{\norm{e^{-2\psi}\Gamma(e^\psi,e^\psi)}_\infty, \norm{e^{2\psi}\Gamma(e^{-\psi},e^{-\psi})}_\infty}\displaybreak[0],\\
    E(t,x,y)&=\sup\set{\abs{\psi(x)-\psi(y)}-t\Lambda(\psi)^2 : \Lambda(\psi)<\infty}.
  \end{align*}
  By Corollary 3.28 of \cite{CaKuStr87} and Lemma \ref{lem:ondiagonalupperboundlargescale} we get for $t\in (0,1]$
  \begin{align*}
    p_{Y^{\rho,\lambda}} (t,x,y) &\leq c t^{-d/\alpha}e^{-E(2t,x,y)} \leq c t^{-d/\alpha}e^{-\abs{\psi(x)-\psi(y)}+2t\Lambda(\psi)^2}
  \end{align*}
  for every $\psi$ with $\Lambda(\psi)<\infty$. Fix $x,y$. Take $\psi (\xi)= \langle \frac{x-y}{\abs{x-y}}, \xi\rangle$. Then $\psi(x)-\psi(y) = \abs{x-y}$. Also we have $\abs{\psi(\xi)-\psi(\zeta)}\leq \abs{\xi-\zeta}$. Now $\abs{e^s-1}\leq \lambda e^\lambda \abs{s}$ for $\abs{s}\leq\lambda$. By (A2) we get
  \begin{align*}
    e^{\mp 2\psi(\xi)}\Gamma(e^{\pm\psi},e^{\pm\psi})(\xi) &= e^{\mp 2\psi(\xi)}\sum_{\substack{\zeta\in\scaleZ{\rho}\\\abs{\xi-\zeta}\leq \lambda}} \big(e^{\pm\psi(\zeta)}-e^{\pm\psi(\xi)}\big)^2 C^\rho(\zeta,\xi) \rho^{-d} \displaybreak[0]\\
    &= \sum_{\substack{\zeta\in\scaleZ{\rho}\\\abs{\xi-\zeta}\leq \lambda}} \big(e^{\pm(\psi(\zeta)-\psi(\xi))}-1\big)^2 C^\rho(\zeta,\xi) \rho^{-d}\\
    &\leq\kappa_1\lambda^2e^{2\lambda}\sum_{\substack{\zeta\in\scaleZ{\rho}\\\abs{\xi-\zeta}\leq \lambda}} \abs{\xi-\zeta}^{2-d-\alpha}\rho^{-d}\leq c_4\lambda^{4-\alpha}e^{2\lambda}\displaybreak[0]
  \end{align*}
  where $c_4>0$ is independent of $x$ and $y$.
\end{proof}

\subsection{Tightness}\label{sect:tightness}

We use the upper bounds established in the previous section to estimate exit times. For this we proceed as in \cite{ChenKumagaiHeatKernelStable} or \cite{BaLe02b} and use the truncated process $Y^{\rho,\lambda}$ together with a standard perturbation argument.

\begin{theorem}\label{theo:tightness-prelim}
  For any $a>0$, $b\in (0,1)$ there exists a constant $\gamma= \gamma(a,b,\kappa_1, \kappa_2, N_0, d,\alpha)$ such that for all $R\geq 1$
  \begin{align*}
    \Prob^x (\tau(B^1(x,aR); Y)<\gamma R^\alpha) &\leq b \quad \forall\, x \in \Z^d \quad \text{ and, equivalently, } \\
    \Prob^x (\tau(B^1(x,aR); Y^\rho)<\gamma R^\alpha) &\leq b \quad \forall\, x \in \scaleZ{\rho}, \forall \, \rho \geq 1\,.
  \end{align*}
\end{theorem}
\begin{proof}
  Fix $\lambda\geq \Theta_2$. First, note that there exists by Lemma \ref{lem:offdiagonalboundtruncated} a constant $c_1>0$ such that for all $t\in[1/2, 1]$
  \begin{align*}
    p_{Y^{\rho,\lambda}} (t,x,y) &\leq c_1 e^{-\abs{x-y}}.
  \end{align*}
  Therefore we estimate for $t\in[1/2, 1]$ and $r>0$
  \begin{align}\label{eq:lemproof1}
    \Prob^x \big(\oabs{Y^{\rho,\lambda}_t-x}\geq r\big) &=\sum_{\substack{y\in\scaleZ{\rho}\\\abs{x-y}\geq r}} p_{Y^{\rho,\lambda}} (t,x,y) \rho^{-d}\leq c_1\sum_{\substack{y\in\scaleZ{\rho}\\\abs{x-y}\geq r}} e^{-\abs{x-y}}\rho^{-d}\leq c_2 e^{-r/2}.
  \end{align}

  With the help of \eqref{eq:lemproof1} we can now estimate the probability that $Y^{\rho,\lambda}$ has left $B^\rho(x,r)$ before $t\leq 1/2$. Introduce the optional time $T_r := \inf\bigset{ t: \bigabs{Y^{\rho,\lambda}_t-Y^{\rho,\lambda}_0}\geq r}$, i.e. the first time the process has left the $r$-ball around its starting point. Then for all $t\leq 1/2$ we get by \eqref{eq:lemproof1} and the strong Markov property:
  \begin{align*}
    \Prob^x\big(\sup_{s\leq t}\bigabs{Y^{\rho,\lambda}_s-x}\geq r\big) &= \Prob^x\big(T_r\leq t; \bigabs{Y^{\rho,\lambda}_1-x}>\frac r2\big)+\Prob^x\big(T_r\leq t; \bigabs{Y^{\rho,\lambda}_1-Y^{\rho,\lambda}_0}\leq\frac r2\big)\\
    &\leq \Prob^x\big(\bigabs{Y^{\rho,\lambda}_1-x}>\frac r2\big) + \int_{0}^t\Prob^x\big(\bigabs{Y^{\rho,\lambda}_1-Y^{\rho,\lambda}_s}>\frac r2 ; T_r\in ds\big)\displaybreak[0]\\
    &\leq c_2e^{-r/4} + \int_{0}^t\E^x\Big(\Prob^{Y^{\rho,\lambda}_s}\big(\bigabs{Y^{\rho,\lambda}_{1-s}-Y^{\rho,\lambda}_0}>\frac r2 ; T_r\in ds\big)\Big)\displaybreak[0]\\
    &\leq c_2e^{-r/4}+c_2e^{-r/4}\Prob^x(T_r\leq t)\leq c_3e^{-r/4}\displaybreak[0].
  \end{align*}

  Hence we have shown for all $t\in [0,1/2]$ that
  \begin{align}\label{eq:truncatedtightness}
    \Prob^x\big(\sup_{s\leq t}\bigabs{Y^{\rho,\lambda}_s-x}\geq r\big) \leq c_3e^{-r/4}.
  \end{align}

  We will now pass over to the process $(Y^\rho_t)$ by handling its large jumps as perturbations of $(Y^{\rho,\lambda}_t)$ with standard techniques of perturbation theory for semigroups, see \cite{KatoPerturbation} for the general case and \cite{Le72} for Markov semigroups. Let $\mathscr{A}^{(\rho, \lambda)}$ resp. $\mathscr{L}^{(\rho)}$ be the generator of $(Y^{\rho,\lambda}_t)$ resp. $(Y^\rho_t)$ and $Q_t^{(\rho,\lambda)}$ resp. $P_t^{(\rho)}$ be the corresponding semigroups. Then
  \begin{align*}
    \mathscr{L}^{(\rho)} = \mathscr{A}^{(\rho, \lambda)} + \mathscr{B}^{(\rho, \lambda)}
  \end{align*}
  where
  \begin{align*}
    \mathscr{B}^{(\rho, \lambda)}f(x) = \sum_{\substack{x,y\in \scaleZ{\rho}\\\abs{x-y}>\lambda}} \big(f(y) - f(x)\big) C^\rho(x,y) \rho^{-d}.
  \end{align*}
  This is a bounded operator in $L^\infty (\scaleZ{\rho})$ since for bounded $f$ by (A2)
  \begin{align*}
    \bigabs{\mathscr{B}^{(\rho, \lambda)}f(x) } &\leq \sum_{\substack{x,y\in \scaleZ{\rho}\\\abs{x-y}>\lambda}} \bigabs{f(y) - f(x)} C^\rho(x,y) \rho^{-d}\leq 2\kappa_1\norm{f}_\infty \sum_{\substack{x,y\in \scaleZ{\rho}\\\abs{x-y}>\lambda}} \abs{x-y}^{-d-\alpha}\rho^{-d}\\
    &\leq c_4\norm{f}_\infty.
  \end{align*}
  Define $S_0(t) = Q_t^{(\rho,\lambda)}$ and for $k>0$ inductively
  \begin{align*}
    S_{k+1} (t) = \int_0^t S_{k}(t-s)\mathscr{B}^{(\rho, \lambda)}S_0(s) ds.
  \end{align*}
  Then it is immediate by the contraction property of $Q_t^{(\rho, \lambda)}$ that $S_k(t)$ is again bounded in $L^\infty$ with operator norm
  \begin{align*}
    \bignorm{S_k(t)}_{\infty,\infty} &\leq \frac{(c_4)^k t^k}{k!}.
  \end{align*}
  Clearly $\sum_{k=0}^\infty S_k(t)$ is well-defined and equals the perturbed semigroup $P_t^{(\rho)}$. We have $\bignorm{P_t^{(\rho)}-Q_t^{(\rho,\lambda)}}_{\infty,\infty}= \bignorm{\sum_{k=0}^\infty S_k(t)}_{\infty,\infty}\leq c_4 te^{c_4t}$ and can therefore estimate for $t\in[1/2,1]$
  \begin{align*}
    \Prob^x\big(\oabs{Y_t^\rho-x}\geq r\big) &= P^{(\rho)}_t \mathds{1}_{B^\rho(x,r)^c}\leq Q^{(\rho,\lambda)}_t \mathds{1}_{B^\rho(x,r)^c}+ c_4te^{c_4t} = \Prob^x\big(\oabs{Y_t^{\rho,\lambda}-x}\geq r\big) + c_4te^{c_4t} \\
    &\leq c_5 e^{-r/2}+ c_5t.
  \end{align*}

  We now proceed as above. Set $T_r := \inf\bigset{ t: \bigabs{Y^{\rho}_t-Y^{\rho}_0}\geq r}$. Then for all $t\leq 1$ we get again by the strong Markov property:
  \begin{align*}
    \Prob^x\big(\sup_{s\leq t}\bigabs{Y^{\rho}_s-x}\geq r\big) &= \Prob^x\big(T_r\leq t; \bigabs{Y^{\rho}_1-x}>\frac r2\big)+\Prob^x\big(T_r\leq t; \bigabs{Y^{\rho}_1-Y^{\rho}_0}\leq\frac r2\big)\\
    &\leq \Prob^x\big(\bigabs{Y^{\rho}_1-x}>\frac r2\big) + \int_{0}^t\Prob^x\big(\bigabs{Y^{\rho}_1-Y^{\rho}_s}>\frac r2 ; T_r\in ds\big)\displaybreak[0]\\
    &\leq c_5e^{-r/4}  + c_5t +\int_{0}^t\E^x\Big(\Prob^{Y^{\rho}_s}\big(\bigabs{Y^{\rho}_{1-s}-Y^{\rho}_0}>\frac r2 ; T_r\in ds\big)\Big)\displaybreak[0]\\
    &\leq c_5e^{-r/4}  + c_5t +\big(c_5 e^{-r/4}+ c_5t \big)\Prob^x(T_r\leq t)\displaybreak[0]\leq c_6e^{-r/4} + c_6t.
  \end{align*}

  This translates by scaling into the following estimate for the process $Y$:
  \begin{align*}
    \Prob^x \big( \sup_{s\leq \rho^\alpha t}\bigabs{Y_s-x}>\rho r\big)\leq c_6 e^{-r/4} + c_6 t.
  \end{align*}

  Given $a$ and $b$ we choose $r>0$ and $t<1/2$ such that the left hand side is smaller than $b$ and in addition $a/r\geq 1$. Now setting $\rho = aR/r$ proves our claim with $\gamma= a^\alpha t/r^\alpha$.
\end{proof}

For later purposes, choose $\tilde\gamma= \gamma(1, 1/2)$, i.e.
\begin{align}\label{eq:choosegamma}
  \Prob^x (\tau(B^1(x,r); Y)<\tilde\gamma r^\alpha)\leq \frac 12\,.
\end{align}

\section{Hitting time estimates and the regularity of the heat kernel}\label{sect:regularity}

In this section we derive an equicontinuity result for the heat kernels of the processes $Y^\rho$. In our application it is essential that the constants appearing do not depend on the scaling parameter $\rho\geq 1$. Again, our presentation uses results from \cite{BaLe02b} and \cite{ChenKumagaiHeatKernelStable}. Another option would be to adopt methods of \cite{Kom95}.

First observe the following L\'{e}vy system identity,  cf. \cite{ChenKumagaiHeatKernelStable}:
\begin{lemma}\label{lem:levysystemidentity}
  Let $f\colon\R^+\times\scaleZ{\rho}\times\scaleZ{\rho} \to \R^+$ be a bounded measurable function vanishing on the diagonal, i.e. $f(t,x,x) =0$ for all $x\in\scaleZ{\rho}$. Then for all $x\in\scaleZ{\rho}$ and predictable stopping times $T$ we have
  \begin{align*}
    \E^x\Big[\sum_{s\leq T} f(s, Y^\rho_{s-}, Y^\rho_{s})\Big] &= \E^x\Big[\int_0^T \Big(\sum _{y\in \scaleZ{\rho} }f (s,Y_s^\rho, y) C^\rho(Y_s^\rho, y)\rho^{-d}\Big) ds\Big].
  \end{align*}
\end{lemma}

Let $W^\rho=(W_t^\rho)_t$ be the space-time process on $\R^+ \times \scaleZ{\rho}$ associated to $Y^{\rho}$, i.e. $W_t = (U_t, Y_t^\rho)$ where $U_t = U_0 +t$ is a deterministic process. We call a measurable function $u\colon\R^+\times \scaleZ{\rho} \to \R$ space-time harmonic or caloric on an open set $\Omega \subset \R^+\times \scaleZ{\rho}$ if for all open relative compact sets $\Omega^\prime\subset\Omega$, $(t,x)\in\Omega^\prime$
\begin{align*}
  u(t,x) = \E^{(t,x)} \big(u(W_{\tau(\Omega; W_s^\rho)}^\rho)\big).
\end{align*}

Important examples for space-time harmonic functions are given by the heat kernel of $Y^\rho$:
\begin{lemma}\label{lem:heatkernelparabolic}
  Let $t_0>0$, $y\in\scaleZ{\rho}$. Then the function $u(t,x) = p_{Y^\rho} (t_0-t,x,y)$ is space-time harmonic in $[0,t_0)\times \scaleZ{\rho}$.
\end{lemma}
The proof is exactly the same as in Lemma 4.5 of \cite{ChenKumagaiHeatKernelStable}.

Define $Q^\rho (t,x,r) := [t, t+\tilde\gamma r^\alpha]\times B^\rho(x,r)$ where $\tilde\gamma$ is chosen such that \eqref{eq:choosegamma} holds. We have the following estimate on the probability of hitting relatively large sets before exiting $Q^\rho (0,x,r)$:
\begin{lemma}\label{lem:hittinglemma1}
  Assume (A1), (A2), (A3) and (A4). Then there exists a constant $c(\kappa_3,d,\alpha)>0$ such that for any $r>\Theta_1\rho^{-1}$, $x\in\scaleZ{\rho}$ and compact set $A\subset Q^\rho(0,x, r)$ with $\leb\otimes \mu^\rho(A) \geq\frac 13 \leb\otimes\mu^\rho \big(Q^\rho(0,x,r)\big)$
  \begin{align*}
    \Prob^{(0,x)}\big(\sigma(A; W^\rho_t) < \tau(Q^\rho(0,x,r);W^\rho_t)\big)\geq c.
  \end{align*}
\end{lemma}

{\bf Remark:} More general, with the proof below and $\Theta_1$, $\gamma$ as in Lemma \ref{lem:necessaryA2A3} we can get lower bounds on the probability of hitting sets $A\subset Q^\rho(0,x,r)$ with $\leb\otimes \mu^\rho(A) \geq2(1-\gamma)\leb\otimes\mu \big(Q^\rho(0,x,r)\big)$ before exiting $Q^\rho(0,x,r)$ for $r\geq \Theta_1\rho^{-1}$ only with assumptions (A1)--(A3). Unluckily, for technical reasons we need (A4) to prove our equicontinuity result.

\begin{proof}
  Set $\tau_r = \tau(Q^\rho(0,x,r); W^\rho_t)$, $\sigma_A=\sigma(A; W^\rho_t)$ and $T = \sigma_A \wedge \tau_r$. Without loss of generality we may assume $\Prob^{(0,x)}(\sigma_{A} < \tau_r)\leq \frac 14$. For each $s\in [0,\infty)$ let $A_s$ denote the projection of $A$ on $\set{s}\times \scaleZ{\rho}$ and let for $y\in \scaleZ{\rho}$
  \begin{align*}
    N(y) = \set{z\in\scaleZ{\rho} : C(y,z)\leq \kappa_3 \abs{y-z}^{-d-\alpha}}.
  \end{align*}

  Choosing $f(s,\xi,\zeta) = \mathds{1}_{Q^\rho (0,x,r)\times A_s\setminus\{(y,y): y\in\scaleZ{\rho}\}} (\xi,\zeta)$ in the L\'{e}vy system formula implies
  \begin{align*}
    \Prob^{(0,x)} \big(\sigma_{A}<\tau_r\big) &\geq \Prob^{(0,x)} \big(\sigma_{A}<\tau_r; Y_{\sigma_A-}\neq Y_{\sigma_A}\big) = \E^{(0,x)}\Big[\sum_{s\leq T}\mathds{1}_{Q^\rho (0,x,r)\times A_s} (Y_{s-},Y_s)\mathds{1}_{\{Y_{s-}\neq Y_s\}}\Big]\\
    &= \E^{(0,x)}\Big[\int_0^T \Big(\sum_{z\in A_s} C^\rho(Y_s^\rho, z)\rho^{-d}\Big)\,ds\Big]\\
    &\geq \E^{(0,x)}\Big[\int_0^T \Big(\sum_{\substack{z\in A_s\setminus N(Y_s^\rho} )} C^\rho(Y_s^\rho, z)\rho^{-d}\Big)\,ds\Big]\\
    &\geq \kappa_3 \E^{(0,x)}\Big[\int_0^T \Big(\sum_{z\in A_s\setminus N(Y_s^\rho)} \abs{Y_s^\rho- z}^{-d-\alpha}\rho^{-d}\Big)\,ds\Big]\displaybreak[0]\\
    &\geq 2^{-d-\alpha}\kappa_3 \E^{(0,x)}\Big[\int_0^T \Big(\sum_{z\in A_s\setminus N(Y_s^\rho)} r^{-d-\alpha}\rho^{-d}\Big)\,ds\Big]\displaybreak[0]\\
    &\geq 2^{-d-\alpha}\kappa_3 r^{-d-\alpha}\E^{(0,x)}\Big[\int_0^T \mu^\rho(A_s\setminus N(Y_s^\rho)) \,ds\Big]\displaybreak[0]\\
    &\geq 2^{-d-\alpha}\kappa_3 r^{-d-\alpha}\E^{(0,x)}\Big[\int_0^{\frac56\tilde\gamma r^\alpha} \mu^\rho(A_s\setminus N(Y_s^\rho)) \,ds\,;\, \sigma_A\wedge\tau_r\geq \tfrac56\tilde\gamma r^\alpha\Big ]\displaybreak[0]\\
    &\geq \frac1{36} \cdot2^{-d-\alpha}\kappa_3 r^{-d-\alpha}(\leb\otimes \mu)(Q^\rho(0,x, r)) \Prob^{(0,x)}\big(\sigma_A\wedge\tau_r\geq \tfrac 56\tilde\gamma r^\alpha\big)\displaybreak[0]\\
    &\geq c(\kappa_3,d,\alpha)\Prob^{(0,x)}\big(\sigma_A\wedge\tau_r\geq \tfrac 56\tilde\gamma r^\alpha\big)\displaybreak[0].
  \end{align*}
  Here we have used $(\leb\otimes \mu)(Q^\rho(0,x, r)) \asymp r^{d+\alpha}$. For the second last step note that for every path with $\sigma_A\wedge\tau_r\geq \tfrac 56\tilde\gamma r^\alpha$
  \begin{align*}
    \frac 13\leb\otimes&\mu^\rho(Q^\rho(0,x, r))\leq \leb\otimes\mu^\rho(A) \\
    &\leq \int_0^{\frac 56\tilde\gamma r^\alpha} \mu^\rho(A_s\setminus N(Y^\rho_s)) ds +  \int_0^{\frac 56\tilde\gamma r^\alpha} \mu^\rho(A_s\cap N(Y^\rho_s))ds\, + \frac 16\leb\otimes\mu^\rho(Q^\rho(0,x, r)).
  \end{align*}
  Now we get by (A4) and because of $r\geq \Theta_3\rho^{-d}$
  \begin{align*}
    \int_0^{\frac 56\tilde\gamma r^\alpha} \mu^\rho(A_s\cap N(Y_s^\rho)) ds &\leq \int_0^{\frac 56\tilde\gamma r^\alpha} \mu(B^\rho(Y^\rho_s,r)\cap N(Y_s^\rho))ds \leq\frac 16\int_0^{\frac56\tilde\gamma r^\alpha} \mu^\rho \big(B^\rho(Y_s^\rho,r)\big)ds\\
    &= \frac 5{36} \tilde\gamma r^\alpha \mu^\rho \big(B^\rho(x,r)\big) = \frac 5{36}\leb\otimes \mu^\rho(Q^\rho(0,x, r)).
  \end{align*}
  Hence
  \begin{align*}
    \int_0^{\frac56\tilde\gamma r^\alpha} \mu^\rho(A_s\setminus N(Y_s^\rho)) ds \geq \frac 1{36}\leb\otimes \mu^\rho(Q^\rho(0,x, r)).
  \end{align*}

  By our choice of $\tilde\gamma$ we obtain
  \begin{align*}
    \Prob^{(0,x)} (\tau_r<\tfrac 56\tilde\gamma r^\alpha)\leq \Prob^x\big(\tau(B^\rho(x,r);Y^\rho)\leq \tilde\gamma r^\alpha\big) \leq \frac 12.
  \end{align*}

  Finally we estimate
  \begin{align*}
    \Prob^{(0,x)}(\sigma_A\wedge \tau_r\geq \tfrac 56\tilde\gamma r^\alpha)&=1-\Prob^{(0,x)}(\sigma_A\wedge \tau_r  <\tfrac 56\tilde\gamma r^\alpha)\\
    &\geq 1- \Prob^{(0,x)} ( \sigma_A<\tau_r) - \Prob^{(0,x)}(\tau_r< \tfrac 56\tilde\gamma r^\alpha)\geq \frac 14.
  \end{align*}
\end{proof}

We also need the following upper bound on the probability of exiting a ball by large jumps.
\begin{lemma}\label{lem:jumpexistestimate}
  Assume (A1), (A2) and (A3). Let $\Theta_1\geq 1$ be the constant of Lemma \ref{lem:necessaryA2A3}. Then there exists a constant $c(\kappa_1,\kappa_2,N_0,d,\alpha)>0$ such that for all $\rho\geq 1$, $s>2r>\Theta_1 \rho^{-1}$, $(t,x)\in [0,\infty)\times\scaleZ{\rho}$
  \begin{align*}
    \Prob^{(t,x)} \big(W^\rho_{\tau(Q^\rho(t,x, r); W^\rho)}\notin Q(t,x, s)\big) \leq c \frac {r^\alpha}{s^\alpha}.
  \end{align*}
\end{lemma}
\begin{proof}
  Let $\tau = \tau (B^\rho(x,r); Y^\rho)$. Observe that $\abs{\xi-\zeta}\geq \frac 12 \abs{x-\zeta}$ for $\xi\in B^\rho(x,r)$ and  $\zeta\in B^\rho(x,s)^c$. Since the space-time process moves continuously in time, it can only exit $Q^\rho(t,x,s)$ and $Q^\rho(t,x,r)$ simultaneously by jumping in space. Using this fact together with the L\'{e}vy system identity for $(Y^\rho_t)$ and (A2) one obtains
  \begin{align*}
    \Prob^{(t,x)}\big(&W^\rho_{\tau(Q^\rho(t,x, r); W^\rho)}\notin Q^\rho(t,x, s)\big)=\Prob^x\big(Y^\rho_{\tau}\notin B^\rho(x,s);\tau\leq \tilde\gamma r^\alpha\big)\leq \Prob^x\big(Y^\rho_{\tau}\notin B^\rho(x,s)\big)\\
    &= \E^x \Big[\sum_{t\leq \tau} \mathds{1}_{\{Y^\rho_{t-}\in B^\rho(x,r), Y^\rho_t\notin B^\rho(x,s)\}}\Big]= \E^x \Big[\int_0^\tau\Big(\sum_{\abs{y-x}\geq s}C^\rho(Y^\rho_t, y)\rho^{-d}\Big)dt \Big]\\
    &\leq \E^x \Big[\int_0^\tau\sum_{\abs{y-x}\geq s}\kappa_1\abs{Y^\rho_t - y}^{-d-\alpha}\rho^{-d}dt \Big]\leq \kappa_1 2^{d+\alpha}\E^x \Big[\tau\sum_{\abs{y-x}\geq s}\abs{x- y}^{-d-\alpha}\rho^{-d} \Big]\\
    &\leq c_1  s^{-\alpha}\E^x (\tau).
  \end{align*}

  Therefore it remains to estimate $\E^x(\tau)$. Let $\kappa_3$ and $\gamma$ be the constants of Lemma \ref{lem:necessaryA2A3}. For each $\xi \in \scaleZ{\rho}$ define $A(\xi) = \bigset{\zeta\in\scaleZ{\rho}: C^\rho(\xi,\zeta)\geq \kappa_3\abs{\xi-\zeta}^{-d-\alpha}}$. We get by applying the L\'{e}vy system identity in the above fashion for $r\geq \Theta_1 \rho^{-1}$
  \begin{align*}
    1&=\Prob^x\big(Y^\rho_\tau \notin B^\rho(x, r)\big)    =\E^x \Big[\sum_{t\leq \tau} \mathds{1}_{\{Y^\rho_{t-}\in B(x,r), Y^\rho_t\notin B(x,r)\}}\Big]\\
     &= \E^x \Big[ \int_0^\tau\Big(\sum_{\abs{y-x}\geq r}C^\rho(Y^\rho_t, y)\rho^{-d}\Big)dt \Big]\geq\E^x \Big[ \int_0^\tau\Big(\sum_{\oabs{y-Y_t^\rho}\geq 2r}C^\rho(Y^\rho_t, y)\rho^{-d}\Big)dt \Big]\\
     &\geq \E^x \Big[ \int_0^\tau\Big(\sum_{\substack{\oabs{y-Y_t^\rho}\geq 2r\\y\in A (Y^\rho_t)}}C^\rho(Y^\rho_t, y)\rho^{-d}\Big)dt \Big]\geq \kappa_3\E^x \Big[ \int_0^\tau\Big(\sum_{\substack{\oabs{y-Y_t^\rho}\geq 2r\\y\in A (Y^\rho_t)}}\oabs{Y_t^\rho -y}^{-d-\alpha} \rho^{-d}\Big)dt \Big]\\
     &\geq \kappa_3\E^x \Big[ \int_0^\tau\Big(\sum_{\substack{2r\leq \oabs{y-Y_t^\rho}< 6r\gamma^{-1/d}\\y\in A (Y^\rho_t)}}\oabs{Y_t^\rho -y}^{-d-\alpha} \rho^{-d}\Big)dt \Big]\\
     &\geq \kappa_3 7^{-d-\alpha}\gamma^{1+d/\alpha}\E^x\Big[\int_0^\tau \rho^{-d} \mu \big(A (Y^\rho_t)\cap B^\rho (Y^\rho_t, 6r\gamma^{-1/d})\setminus B^\rho (Y^\rho_t, 2r)\big) dt\Big]\\
     &\geq c_1(\kappa_3,d,\alpha)\E^x\Big[\int_0^\tau \rho^{-d} \mu \big(B^\rho (Y^\rho_t, 2r)\big) dt\Big]\geq c_2(\kappa_3,d,\alpha)r^{-\alpha}  \E^x(\tau).
  \end{align*}
  Here the second last inequality is due to Lemma \ref{lem:necessaryA2A3} since there are at least $2\mu\big(B^\rho (Y^\rho_t, 2r)\big)$ elements of $A(Y_t^\rho)$ in $B^\rho (Y^\rho_t, 6r\gamma^{-1/d})$.
\end{proof}

\begin{proposition}\label{prop:regularityparabolic}
  Assume (A1)-(A4). Then there exist constants $c>0$ and $\beta \in (0,1)$ depending only on the constants appearing in (A1)-(A4) such that for all $R>\Theta_1\rho^{-1}$, $q\colon [0, \tilde\gamma 16R]\times \scaleZ{\rho} \to \R$ bounded and space-time harmonic in $Q^\rho(0, x_0, 16R)$ the following a-priori continuity estimate holds:
  \begin{align*}
    \oabs{q(s,x)-q(t,y)}\leq c\norm{q}_{\infty} R^{-\beta}\big(\abs{s-t}^{1/\alpha} + \abs{x-y}\big)^\beta
  \end{align*}
  for all $(s,x), (t,y)\in Q^\rho(0, x_0,R)$ with $\abs{x-y}\geq \Theta_1\rho^{-1}$.

  Moreover, if $\abs{x-y}\leq \Theta_1\rho^{-1}$ we have
  \begin{align*}
     \oabs{q(s,x)-q(t,y)}\leq c\norm{q}_{\infty} R^{-\beta}\big(\abs{s-t}^{1/\alpha} + \Theta_1\rho^{-1}\big)^\beta.
  \end{align*}
\end{proposition}
\begin{proof}
  The proof can be found in \cite{BaKu06} or \cite{ChenKumagaiHeatKernelStable} for our case and in \cite{BaKa05} or \cite{HuKa06a} for the ``elliptic'' case.

  We may assume $\norm{q}_\infty = 1/2$. Fix $(t,x)\in Q^\rho(0, x_0, R)$. Let $\xi\in (0,1/2)$, $\eta>0$, $Q_k = Q^\rho(t,x_0, \xi^kR)$, $\tau_k = \tau(Q_k; Z^\rho)$ and define for $k\in\N$
  \begin{align*}
    m_k &=\inf_{z\in Q_k} q(z), &M_k &=\sup_{z\in Q_k} q(z).
  \end{align*}

  We show that it is possible to choose constants $\xi, \zeta$ independent of $R>0$, $x,x_0\in Q^\rho( 0, x_0, R)$ and $q$ such that
  \begin{align}\label{eq:HoelderproofCentralAim}
    M_k- m_k \leq \zeta^k
  \end{align}
  for all $k\geq 0$ with $\xi^kR\geq \Theta_1\rho^{-1}$ where $\Theta_1$ the constant in (A4). The restriction that $x$ and $y$ cannot be arbitrarily close is natural since hitting time estimates in the sense of Lemma \ref{lem:hittinglemma1} may not hold for small $r$. Just consider the conductivities $C(x,y)= \abs{x-y}^{-d-\alpha}\mathds{1}_{\set{\abs{x-y}>R_0}}$. Then the process $(W^t)$ has only jumps with length bigger than $R_0$ and therefore can't hit sets $A\subset B(x_0, r)$ for $r<R_0$ starting in $x_0$ unless $x_0\in A$.

  Trivially, \eqref{eq:HoelderproofCentralAim} holds for $k= 0$. Now assume that this equation holds already for all $i\leq k$ while still $\xi^{k+1}R\geq \Theta_1\rho^{-1}$. We set
  \begin{align*}
    A_k  = \bigset{z\in Q_{k+1} : q(z) \leq \frac{M_k+m_k}{2}}.
  \end{align*}
  Without loss of generality we might assume $\leb\otimes \mu(A_k)/\leb\otimes\mu({Q_{k+1}})\geq 1/2$. Else we just look at $1/2-q$ instead of $q$. We choose a compact set $A^\prime_k \subset A_k$ with $\leb\otimes \mu(A_k^\prime)/\leb\otimes\mu({Q_{k+1}})\geq 1/3$ and define $T_k = T(A^\prime_k; W^\rho)$.   Observe that $u$ is harmonic in $Q_k\subset Q^\rho(t,x, R)\subset Q^\rho(0,x_0, 16R)$. Therefore we get for arbitrary $z_1, z_2\in Q_{k+1}$
  \begin{align*}
    q(z_1)-q(z_2) &=  \E^{z_1} \big[q(W_{T_k\wedge\tau_{k+1}}^\rho)\big]- q(z_2)\\
    &= \E^{z_1} \big[q(W_{T_k}^\rho)- q(z_2); T_k<\tau_{k+1}\big]  \\
    &\quad+ \E^{z_1} \big[q(W_{\tau_{k+1}})- q(z_2); T_k>\tau_{k+1}\text{ and } W_{\tau_{k+1}}^\rho\in Q_k\big] \\
    &\quad+ \sum_{i=1}^k \E^{z_1}\big[q(W_{\tau_{k+1}})- q(z_2); T_k>\tau_{k+1}\text{ and } W_{\tau_{k+1}}^\rho\in Q_{k-i}\setminus Q_{k-i+1}\big]\\
    &\quad + \E^{z_1} \big[q(W_{\tau_{k+1}})- q(z_2); T_k>\tau_{k+1}\text{ and } W_{\tau_{k+1}}^\rho\notin Q_0\big].
  \end{align*}

  Here, the first term can be estimated by $\frac12(M_k-m_k)\Prob^{z_1} (T_k> \tau_{k+1})$ while the second term is bounded from above by $(M_k-m_k)\Prob^{z_1} (T_k> \tau_{k+1})= (M_k-m_k)(1-\Prob^{z_1} (T_k< \tau_{k+1}))$. Moreover Lemma \ref{lem:hittinglemma1} implies the existence of $c_1>0$ such that $\Prob^{z_1} (T_k<\tau_{k+1})\geq c_1$. For the remaining terms note that by Lemma \ref{lem:jumpexistestimate} there exists $c_2>0$ with
  \begin{align*}
    \Prob^{z_1} (W^\rho_{\tau_{k+1}}\notin Q_i) \leq c_2\xi^{(k+1-i)\alpha}
  \end{align*}

  Summing up the terms above we this get for $\xi^\alpha\leq\zeta/2$
  \begin{align*}
    q(z_1)-q(z_2) &\leq (M_k-m_k)\big(1-\frac 12\Prob^{z_1}(T_k<\tau_{k+1})\big) + \sum_{i=1}^{k+1} (M_{k-i}-m_{k-i})\Prob^{z_1} (Z^\rho_{\tau_{k+1}}\notin Q_{k-i+1}) \displaybreak[0]\\
    &\leq (1-\frac 14 c_1)\zeta^k -\frac 14 c_1\zeta^k + c_2\sum_{i=1}^{k+1}\zeta^{k-i} \xi^{\alpha i}\displaybreak[0]\\
    &\leq (1-\frac 14 c_1)\zeta^k -\frac 14 c_1\zeta^k + c_2\zeta^k\sum_{i=1}^{\infty} \big(\frac{\xi^\alpha}\zeta\big)^{i}
    =(1-\frac 14 c_1)\zeta^k -\frac 14 c_1\zeta^k + 2c_2\xi\zeta^k.\displaybreak[0]
  \end{align*}

  Estimate \eqref{eq:HoelderproofCentralAim} follows from the equation above by taking $\zeta = 1-c_1/4$ and $\xi = \frac 12 \wedge \big(\frac\zeta2\big)^{1/\alpha}\wedge \frac {c_1}{8c_2}$.

  To derive H\"{o}lder continuity let $z_i= (s_i, x_i)\in Q^\rho(0, x_0, R)$, $i=1, 2$ with $z_1\neq z_2$ and $s_1\leq s_2$. Assume $\abs{x_1-x_2}\geq \Theta_1\rho^{-1}$ and take $k$ maximal such that $z_2\in B^\rho(z_1, R\xi^k)$. Then
  \begin{align*}
    \abs{s_1-s_2}^{1/\alpha}+\abs{x_1-x_2}&\leq (\tilde\gamma^{1/\alpha} +1) R\xi^{k}, &\xi^kR&\geq \Theta_1\rho^{-1}, &\abs{q(z_1)-q(z_2)}&\leq \zeta^k.
  \end{align*}

  Thus by optimality $k$ is the smallest integer such that \[k\geq (\log \xi)^{-1}\big(\log (\abs{s_1-s_2}^{1/\alpha}+\abs{x_1-x_2})-\log(\tilde\gamma^{1/\alpha}+1) R\big)-1,\] and we get
  \begin{align*}
    \abs{q(z_1)-q(z_2)} &\leq \zeta^{-1}(\tilde\gamma^{1/\alpha} +1)^{-\log\zeta/\log\xi}R^{-\log\zeta/\log\xi}(\abs{s_1-s_2}^{1/\alpha} + \abs{x_1-x_2})^{\log \zeta/\log\xi},
  \end{align*}
  i.e. the proposition holds with $\beta = \log \zeta/\log\xi \in (0,1)$ and $c= \zeta^{-1}(\tilde\gamma^{1/\alpha} +1)^{-\beta}$.
  If on the other hand $\abs{x_1-x_2}\leq \Theta_1(\tilde\gamma^{1/\alpha}+1)\rho^{-1}$ we take $k$ maximal with $R\xi^k\geq \Theta_1\rho^{-1}$ and $z_2\in B^\rho(z_1, R\xi^k)$. Then in particular $\abs{s_1-s_2}^{1/\alpha}\leq \tilde\gamma^{1/\alpha} R\xi^{k}$, and we get for $k$ the inequalities
  \begin{align*}
     k&\leq (\log \xi)^{-1} \big(\log(\Theta_1\rho^{-1})-\log R\big),\\
     k&\leq (\log \xi)^{-1} \big(\log(\abs{s_1-s_2}^{1/\alpha})-\log (\tilde\gamma^{1/\alpha} R)\big)
  \end{align*}
  Combining this we get with $\beta$ as above
  \begin{align*}
    \abs{q(z_1)-q(z_2)} &\leq \zeta^{-1}(\tilde\gamma^{-\beta/\alpha} +1)R^{-\beta}\big(\abs{s_1-s_2}^{\beta/\alpha} + (\Theta_1\rho^{-1})^{\beta}\big)\\
    &\leq 2\zeta^{-1}(\tilde\gamma^{-\beta/\alpha} +1)R^{-\beta}\big(\abs{s_1-s_2}^{1/\alpha} + \Theta_1\rho^{-1}\big)^\beta
  \end{align*}
\end{proof}

In particular, this implies regularity of the heat kernels
\begin{theorem}\label{theo:regularityheatkernels}
  There exist constants $c>0$ and $\beta\in (0,1)$ such that for all $t_0\in(0,\infty)$, $x_1, x_2, y\in \scaleZ{\rho}$ and $s_1, s_2\geq [t_0,\infty)$
  For $y\in\scaleZ{\rho}$ we have
  \begin{align*}
    &\abs{p_{Y^\rho}(s_1, x_1, y)-p_{Y^\rho}(s_2, x_2, y)}\\
    &\quad\quad\quad\leq c t_0^{-(d+\beta)/\alpha} \big(\abs{s_1-s_2}^{1/\alpha}+\abs{x_1-x_2}\vee(\Theta_1\rho^{-1})\big)^{\beta}.
  \end{align*}

  More general, for arbitrary $y_1,y_2\in\scaleZ{\rho}$
  \begin{align*}
    &\abs{p_{Y^\rho}(s_1, x_1, y_1)-p_{Y^\rho}(s_2, x_2, y_2)}\\
    &\quad\quad\quad\leq c t_0^{-(d+\beta)/\alpha} \big(\abs{s_1-s_2}^{1/\alpha}+\abs{x_1-x_2}\vee(\Theta_1\rho^{-1})+\abs{y_1-y_2}\vee (\Theta_1\rho^{-1})\big)^{\beta}.
  \end{align*}
\end{theorem}
\begin{proof}
  Fix $t_0>0$. For an arbitrary $T_0\geq t_0$ the function $q(t,x) = p(T_0-t, x,y)$ is space-time harmonic on $[0, T_0/2)\times \scaleZ{\rho}$ by Lemma \ref{lem:heatkernelparabolic} as well as bounded by $c_1t_0^{-d/\alpha}$ by Lemma \ref{lem:ondiagonalupperboundallscales}. Now take $R=(T_0/32\tilde\gamma)^{1/\alpha}$ and $s_1,s_2\in[0,\tilde\gamma R^\alpha)$. Assume first $\Theta_1\rho^{-1}\leq\oabs{x_1-x_2}$. If $\oabs{x_1-x_2}\leq R$ we get by Proposition \ref{prop:regularityparabolic}
  \begin{align*}
     \oabs{p_{Y^\rho}(T-s_1, x_1, y)-p_{Y^\rho}(T-s_2, x_2, y)}&\leq c_2 (T_0/32\tilde\gamma)^{-\beta/\alpha}c_1 t_0^{-d/\alpha}\big(\abs{s_1-s_2}^{1/\alpha}+\abs{x_1-x_2}\big)^\beta\\
     &\leq c_3 t_0^{-(d+\beta)/\alpha} \big(\abs{s_1-s_2}^{1/\alpha}+\abs{x_1-x_2}\big)^\beta.
  \end{align*}

  In the other case $\oabs{x_1-x_2}> R$ we have $(\abs{s_1-s_2}^{1/\alpha}+\abs{x_1-x_2})^\beta\geq c_4 t_0^{\beta/\alpha}$ and hence
  \begin{align*}
    \oabs{p_{Y^\rho}(T-s_1, x_1, y)-p_{Y^\rho}(T-s_2, x_2, y)}&\leq 2c_1 t_0^{-d/\alpha}\leq2c_1 t_0^{-(d+\beta)/\alpha}t_0^{\beta/\alpha}\\
     &\leq c_5 t_0^{-(d+\beta)/\alpha} \big(\abs{s_1-s_2}^{1/\alpha}+\abs{x_1-x_2}\big)^\beta.
  \end{align*}
  In the same fashion we can deal with the case $\Theta_1\rho^{-1}\geq\oabs{x_1-x_2}$. The other a-priori estimate asserted in the theorem now follows by symmetry of the heat kernels.
\end{proof}

\section{The central limit theorem}\label{sect:centrallimittheorem}

The aim of this section is to provide a proof of Theorem \ref{theo:meta-one} and Theorem \ref{theo:meta-two}.  Although the limit process $\sX$ is a jump process the idea of the proof is very similar to the one in \cite {StrZheng97} and \cite{BaKu06}. For $n \in \N$ let $Y^n=(Y^n_t)_t$ be the time-continuous processes defined in the introduction and explained in section \ref{sect:definitions}. Denote the corresponding semigroup by $(P^{(n)}_t)_t$ and its kernel by $p^{(n)}(t,x,y)$. Recall that the Dirichlet form corresponding to $Y^n$ is given by
\begin{align*}
  \mathcal{E}^{(n)} (f,f) =  \sum_{x,y\in \scaleZ{n}} \big(f(y) - f(x)\big)^2 C^n(x,y) n^{-2d}.
\end{align*}
We denote the restriction of functions on $\R^d$ to $\scaleZ{n}$ by $R_n$. We also need to extend functions on the grid to continuous functions on $\R^d$ in a (for our purpose) reasonable way: For $x\in \scaleZ{n}$ set $\mathfrak{Q}_n(x):= \prod [x_i, x_i+1/n)$ and $\mathfrak{Q}_n=\bigcup\limits_{x\in\scaleZ{n}} \{ \mathfrak{Q}_n(x)\}$. The sets of $\mathfrak{Q}_n$ form a partition of $\R^d$. Recall that, for a point $x \in \R^d$ we denote by $[x]_n$ the element of $\scaleZ{n}$ satisfying $([x]_{n})_i=n^{-1}\lceil nx_i \rceil$ for all $i=1,\ldots,d$. For $f\colon \scaleZ{n}\to \R$ we denote by $E_n f$ a Lipschitz-continuous function $E_n f\colon \R^d \to \R$  satisfying:
\begin{align*}
&a) \quad (E_n f)(x) = f(x) \quad \forall x\in \scaleZ{n} \,, \\
&b) \quad \underline{f}(x) \leq E_n f(x)\leq \overline{f}(x) \; \forall\, x \in \R^d\,, \text{ where } \underline{f}(x) = \min\limits_{\overline{\mathfrak{Q}_n(x)}\cap \scaleZ{n}} f;\overline{f}(x) = \max\limits_{\overline{\mathfrak{Q}_n(x)}\cap \scaleZ{n}} f \,, \\
&c) \quad \bigabs{\nabla E_n f(x)} \leq c n \max\big\{|f(\xi)-f(\eta)|: \xi, \eta \in \overline{\mathfrak{Q}_n(x)} \cap \scaleZ{n} \big\} \;  \forall \, n \in \N, \forall\, x \in \R^d \,.
\end{align*}
The precise choice of the function $E_n f$ is not important for our approach as long as $E_n$ is a linear operator.

Let us emphasize that in the following result we adopt the notion $\Prob^x$ for the probability of a Markov process starting in $x$. Any starting point, together with a stochastic process with c\`{a}dl\`{a}g paths, corresponds to a probability measure on $D([0,\infty), \R^d)$. For a Markov process $X=(X_t)$ starting in $x$ we refer to this probability measure as ''the law of $X$ under $\Prob^x\,$''.

The following theorem is a reformulation of Theorem \ref{theo:meta-one}.

\begin{theorem}\label{theo:tightness}
Let $(C^n)_n$ be a sequence of conductivity functions satisfying (A1) through (A4) and $(x_n)_n$ a sequence of points $x_n \in \scaleZ{n}$ with $x_n \to x \in \R^d$ for $n \to \infty$. Then each subsequence $(n^\prime)$ of $(n)$ has a subsequence $(n^{\prime\prime})$ with the following properties:
  \begin{enumerate}
    \item For any  $f\in C_c (\R^d)$ the continuous functions $\big(E_{n^{\prime\prime}} P^{(n^{\prime\prime})}_t R_{n^{\prime\prime}} f\big)$ converge uniformly on compact sets for $n^{\prime\prime}\to \infty$. The limit defines a family of linear operators $(P_t)_{t>0}$ which extends to the semigroup on $C(\R^d)$ of a symmetric strong Markov process $\sX$.
    \item For any $t_0>0$ the laws of $(Y^{n^{\prime\prime}}_t)_{t\in [0,t_0]}$ under $\Prob^{x_{n^{\prime\prime}}}$ converge weakly to the law of $(\sX_t)_{t\in [0,t_0]}$.
  \end{enumerate}
\end{theorem}

Once these assertions are proved it remains to show that $\sX$ does not depend on the choice of $(n')$. The proof of Theorem \ref{theo:tightness} makes use of the following sufficient condition for tightness:

\begin{theorem}[{\cite{Al78}}]\label{thm:Aldous}
  Let $(Y^n)$ be a sequence of stochastic processes with c\`{a}dl\`{a}g paths, $t_0>0$. Assume:
  \begin{enumerate}
    \item For all sequences $\sigma_n$ of random variables with values in $[0,t_0]$ such that $\sigma_n$ is a stopping time with respect to the natural filtration of $Y^n$, sequences $\delta_n\geq 0$ with $\lim\limits_{n\to \infty} \delta_n = 0$ and $\eta>0$ it holds that
    \begin{align}\label{eq:critAldous}
      \lim_{n\to\infty}\Prob\Big(\bigabs{Y^n_{\sigma_n+\delta_n}- Y^n_{\sigma_n}}>\eta\Big) = 0\,.
    \end{align}
    \item Either $(Y^n_0)$ and $\max_{t\in[0,t_0]} \bigabs{Y^n_t-Y^n_{t-}}$ are tight or $Y^n_t$ is tight for every $t\in [0,t_0]$.
  \end{enumerate}
  Then the laws of $(Y^n)_n$ are tight in $D([0,t_0], \R^d)$.
\end{theorem}

\begin{proofo}{Proof of Theorem \ref{theo:tightness} and Theorem \ref{theo:meta-one}}
  Let $(n^\prime)$ be a subsequence of $(n)$. Fix countable dense subsets $(s_i)$ of $[0,\infty)$ and $(f_j)$ of $C_c(\R^d)$.

By $Q^{(n)}_t := E_n P^{(n)}_t R_n$ we define a positivity-preserving contraction semigroup $(Q^{(n)}_t)_t$ on the Banach space $C(\R^d)$. Now,  Theorem  \ref{theo:regularityheatkernels} yields that for all $i$, $j$ the family of functions $(Q^{(n)}_{s_i} f_j)_{ n\in\N}$ is equicontinuous. In fact, we have for $x,y \in \scaleZ{n}$ with $\oabs{x-y}\geq n^{-1} \Theta_1$
  \begin{align*}
    \bigabs{P^{(n)}_{s_i} R_n(f_j)(x)- P^{(n)}_{s_i}  R_n(f_j)(y)}&\leq \sum_{z\in \scaleZ{n}}\bigabs{p^{(n)}(s_i,x,z)-p^{(n)}(s_i,y,z)}\bigabs{f_j(z)}n^{-d}\\
        &\leq cs_i^{-(d+\beta)/\alpha}\leb(\supp(f_j))\norm{f_j}_\infty \abs{x-y}^{\beta}
  \end{align*}
  where $c>0$, $\Theta_1>0$ and $\beta\in (0,1)$ are independent of $n$ and stem from Proposition \ref{prop:regularityparabolic}. The construction of the extension operators $E_n$ implies for all $x,y\in \R^d$ with $\oabs{x-y}\geq n^{-1}\Theta_3$
  \begin{align}\label{eq:equicontinuitylargescale}
    \bigabs{Q^{(n)}_{s_i} f_j (x) - Q^{(n)}_{s_i} f_j (y)}\leq 2cs_i^{-(d+\beta)/\alpha}\leb(\supp(f_j))\norm{f_j}_\infty \abs{x-y}^{\beta}
  \end{align}
In a similar fashion one establishes for $\oabs{x-y}\leq n^{-1}\Theta_1$
  \begin{align}\label{eq:equicontinuitysmallscale}
    \bigabs{Q^{(n)}_{s_i} f_j (x) - Q^{(n)}_{s_i} f_j (y)}\leq 2c\Theta_1s_i^{-(d+\beta)/\alpha}\leb(\supp(f_j))\norm{f_j}_\infty n^{-\beta}.
  \end{align}
  Furthermore  $(Q^{(n)}_{s_i} f_j)$ is equibounded. The Theorem of Arzela-Ascoli and the passage to a diagonal sequence give us therefore a subsequence $(n^{\prime\prime})$ of $(n^\prime)$ such that for all $i$, $j$ the sequence $Q^{(n^{\prime\prime})}_{s_i} f_j$ converges uniformly on compact sets for $n^{\prime\prime}\to \infty$. Denote this limit function by $P_{s_i} f_j$. We use an $\nicefrac{\varepsilon}{3}$-argument to extend it for all positive times: Let $t\in (0,\infty)$ and take a subsequence $(i^\prime)$ of $(i)$ such that $s_{i^\prime}\to t$ for $i^\prime\to \infty$ and $s_{i^\prime}\geq t/2$. Then
  \begin{align*}
    \bigabs{Q^{(n^{\prime\prime})}_t f_j(x)-Q^{(m^{\prime\prime})}_t f_j(x)}&\leq \bigabs{Q^{(n^{\prime\prime})}_t f_j(x)-Q^{(n^{\prime\prime})}_{s_{i^\prime}} f_j(x)}+\bigabs{Q^{(n^{\prime\prime})}_{s_{i^\prime}} f_j(x)-Q^{(m^{\prime\prime})}_{s_{i^\prime}} f_j(x)}\\
      &\quad\quad\quad+\bigabs{Q^{(m^{\prime\prime})}_t f_j(x)-Q^{(m^{\prime\prime})}_{s_{i^\prime}} f_j(x)}.
  \end{align*}

The second term on the right hand side converges uniformly on compact sets to $0$ for $n^{\prime\prime},m^{\prime\prime}\to\infty$. The other two terms can be handled again by Theorem \ref{theo:regularityheatkernels}:
  \begin{align*}
    \bigabs{P^{(n)}_t R_n(f_j)(x)- P^{(n)}_{s_{i^\prime}} R_n(f_j)(x)}&\leq \sum_{z\in \scaleZ{n}}\bigabs{p^{(n)}(t,x,z)-p^{(n)}(s_i,x,z)}\bigabs{f_j(z)}n^{-d} \\
    &\leq c_1s_i^{-(d+\beta)/\alpha}\leb(\supp(f_j))\norm{f_j}_\infty \big(\oabs{s_i-t}^{1/\alpha}+\Theta_3n^{-1}\big)^\beta.
  \end{align*}
  The right hand side clearly converges to $0$ for $n\to\infty$. Hence the limit $P_{t} f_j$ exists uniformly on compact sets for all $t\in [0,\infty)$. Finally by $\bignorm{Q^{(n^{\prime\prime})}_{t} f_j}\leq \norm{f_j}$ and because $(f_j)$ is dense in $C(\R^d)$ in the topology of uniform convergence on compact sets we have established the desired convergence result for all $f\in C(\R^d)$.

  It follows from the corresponding properties of the $Q^{(n)}_t$ that $P_t$ is a positivity-preserving contraction semigroup on $C(\R^d)$ which is hence associated to a symmetric strong Markov process $\sX$ on $\R^d$.

  Fix $t_0>0$ and $x\in\R^d$. We want to apply Theorem \ref{thm:Aldous}. Take an arbitrary sequence of stopping times $\tau_n \in [0,t_0]$, a sequence $(\delta_n)$ of reals  converging to $0$ and $a>0$. By Theorem \ref{theo:tightness-prelim} for each choice $b\in (0,1)$ there exist a constant $\gamma(a, b)$ with
  \begin{align}\label{eq:exittimeestimate}
    \Prob^{x_n} \big(\tau\big(B(x_n, a); Y^{(n)}\big)\leq\gamma(a,b)\big) \leq b
  \end{align}
  for all $n\in\N$. Therefore, for all $n$ large enough such that $\delta_n\leq \gamma(a, b)$,
  \begin{align*}
    \Prob^{x_n} \big(\bigabs{Y^{(n)}_{\tau_n+\delta_n}-Y^{(n)}_{\tau_n}}>a\big) &= \Prob^{x_n} \big(\bigabs{Y^{(n)}_{\delta_n}-Y^{(n)}_{0}}>a\big)\\
    &\leq \Prob^{x_n} \big(\tau\big(B(x_n, a); Y^{(n)}\big)\leq \delta_n\big) \\
    &\leq \Prob^{x_n} \big(\tau\big(B(x_n, a); Y^{(n)}\big)\leq \gamma (a,b)\big)\leq b
  \end{align*}
  where the first equality follows by the strong Markov property. \eqref{eq:critAldous} follows immediately. Moreover, $x_n\to x$ implies the tightness of the starting distributions while \eqref{eq:exittimeestimate} implies the tightness of $\max_{t\in[0,t_0]} \bigabs{Y^{(n)}_t-Y^{(n)}_{t-}}$, both under $\Prob^{x_n}$. The tightness of the laws of $(Y^{(n)}_t)$ under $\Prob^{x_n}$ follows.

  Finally we prove the asserted weak convergence for $n^{\prime\prime}\to\infty$ by showing that the finite dimensional distributions of the limit probability $\mathds{Q}$ on $D([0,t_0], \R^d)$ of a weakly convergent subsequence $(n^{\prime\prime\prime})$ are independent of the actual subsequence. For $g\in C_c (\R^d)$, $t\in[0,t_0]$
  \begin{align*}
     \int g(\omega_t)d\mathds{Q}(\omega) &= \lim _{n^{\prime\prime\prime}\to\infty} \E^{x_{n^{\prime\prime\prime}}} (R_{n^{\prime\prime\prime}} g)\big(X^{(n^{\prime\prime\prime})}_t\big) = \lim _{n^{\prime\prime\prime}\to\infty} P^{(n^{\prime\prime\prime})}_t (R_{n^{\prime\prime\prime}} g)(x_{n^{\prime\prime\prime}})\\
    &= P_t g(x)\,,
  \end{align*}
  where the last equality follows from the equicontinuity of the family $P^{(n^{\prime\prime\prime})}_tg$. Therefore the one-dimensional distributions are independent of $(n^{\prime\prime\prime})$. More generally let $0\leq s_1<\ldots <s_k\leq t_0$ and $g_1,\ldots, g_k\in C_c(\R^d)$. Then by the time-homogeneity of our Markov chains
  \begin{align*}
    \int g_1(\omega_{s_1})\cdot\ldots\cdot g_k(\omega_{s_k})\,d\mathds{Q}(\omega) &=\lim _{n^{\prime\prime\prime}\to\infty}  \E^{x_{n^{\prime\prime\prime}}} \big(g_1(X_{s_1}^{(n^{\prime\prime\prime})})\cdot\ldots\cdot g_k(X_{s_k})^{(n^{\prime\prime\prime})})\big)  \\
    &=\lim_{n^{\prime\prime\prime}\to\infty} \big(P^{(n^{\prime\prime\prime})}_{s_1}\big(g_1 P^{(n^{\prime\prime\prime})}_{s_2-s_1}(\ldots  P^{(n^{\prime\prime\prime})}_{s_k-s_{k-1}}g_k)\big)\big)(x_{n^{\prime\prime\prime}})\\
    &=\big(P_{s_1}\big(g_1 P_{s_2-s_1}(\ldots  P_{s_k-s_{k-1}}g_k)\big)\big)(x)\,.
  \end{align*}
  Here the last equality is again due to the equicontinuity. Hence the k-dimensional distributions of $\mathds{Q}$ are independent of the choice of the subsequence $(n^{\prime\prime\prime})$ and are determined by the semigroup $(P_t)$. Therefore we have weak convergence along $(n^{\prime\prime})$. In particular, the stochastic process corresponding to $\mathds{Q}$ has the same finite-dimensional distributions as $\sX$ starting in $x$.
\end{proofo}

We proceed to the proof of Theorem \ref{theo:meta-two}.

\begin{proofo}{Proof of Theorem \ref{theo:meta-two}}
Let $(n')$ be any subsequence of $(n)$. Let $\sX$ be a strong Markov process - possibly depending on the choice of $(n')$ - and $(n'')$ be a subsequence of $(n')$ such that the assertions of Theorem \ref{theo:tightness} hold true. We aim to show that $\sX$ does not depend on the choice of $(n')$. It suffices to show that the limiting process $\sX$ corresponds to the Dirichlet form \eqref{eq:defE}. This is the case if
  \begin{align}\label{eq:resolventeqlimit}
    \mathcal{E} (U_\lambda f,g) &= (f,g) - \lambda(U_\lambda f, g)
  \end{align}
  for any $f,g\in C^\infty_c(\R^d)$ where $U_\lambda f(x)= \int_0^\infty e^{-\lambda t} (P_t f)(x) dt$, $\lambda>0$. Note that, at this stage, $U_\lambda$ does depend on the choice of $(n')$. Equality (\ref{eq:resolventeqlimit}) implies
  \[ \mathcal{E} (U_\lambda f,g) + \lambda(U_\lambda f, g) = \mathcal{E} (G_\lambda f,g) + \lambda(G_\lambda f, g) \quad \text{ for any } f,g\in C^\infty_c(\R^d)\,, \]
 where $G_\lambda$ and $(\sE, D(\sE))$ are independent of $(n')$. $G_\lambda$ is then the $L^2$-resolvent of $\sX$ and we are done. Note that Theorem \ref{theo:equivnorms} implies
  \begin{align*}
    D(\mathcal{E}) = H^{\alpha/2}(\R^d).
  \end{align*}

  We prove (\ref{eq:resolventeqlimit}) by approximating each term by its discrete analog. On the discrete level
  \begin{align}\label{eq:resolventapprox}
    \mathcal{E}^{(n)} (U^{(n)}_\lambda R_n(f), R_n(g)) &= (R_n (f),R_n (g)) - \lambda (U^{(n)}_\lambda R_n(f), R_n(g))
  \end{align}
  where $U_\lambda^{(n)} h(x)= \int_0^\infty e^{-\lambda t} \big(P^{(n)}_t h\big)(x)\,dt$, $\lambda>0$, denotes the resolvent of $\big(Z_t^{(n)}\big)$.

  Therefore fix $\lambda>0$, $f,g\in C^\infty_c (\R^d)$ and abbreviate $f_n = R_n(f)$, $g_n = R_n(g)$. Then $f_n, g_n \in L^2(\scaleZ{n},n^{-1}\mu)$ with $\norm{f_n}+\norm{g_n} \leq c$ for all $n$. Recalling the definition of section \ref{sect:assumform} one sees that $\sum_{x\in \scaleZ{n}} f_n(x)\mathds{1}_{\mathfrak{Q}_n}(x)$ converges in $L^2(\R^d)$ to $f$ and $\abs{(f_n, g_n) - (f,g)}$ converges to zero for $n\to\infty$.

  Now by the compactness of the support of $f$ and Theorem \ref{theo:regularityheatkernels} we get equicontinuity for the family $(E_nU_\lambda^{(n)} R_n f)_{n\in\N}$ analogous to \eqref{eq:equicontinuitylargescale} resp. \eqref{eq:equicontinuitysmallscale}. Together with $\bigabs{[x]_n - x}\leq \sqrt{d}n^{-1}$ we get
  \begin{align*}
     \bigabs{E_n U_\lambda^{(n)}f_n (x) - U_\lambda^{(n)}f_n([x]_n)} \leq c n^{-\beta}
  \end{align*}
  for all $x\in\R^d$. In particular, the functions $x\mapsto U_\lambda^{(n')} f_{n'} ([x]_{n'})$ on $\R^d$ converge along the subsequence $(n'')$ uniformly on compact sets to $U_\lambda f$. Taking into account that $g$ is compactly supported we get by dominated convergence
  \begin{align}
    \bigabs{\big(U_\lambda^{(n'')} f_{n''}, g_{n''}\big)- \big(U_\lambda f, g\big)}  \to 0 \quad \text{ for }n'' \to\infty.
  \end{align}

  Therefore the right-hand side of \eqref{eq:resolventapprox} converges against the right-hand side of \eqref{eq:resolventeqlimit} for the subsequence $n'' \to\infty$.

  It remains to determine the limit of the left-hand side of \eqref{eq:resolventapprox} for $n'' \to \infty$. We do this in several steps.

  {\bf Step 1:} $U_\lambda f \in H^{\alpha/2} (\R^d)$.

  This result probably follows from standard arguments of approximation theory. For the sake of completeness we give a detailed proof. First note that ${U^{(n)}_\lambda f_n}$ and $E_n U^{(n)}_\lambda f_n$ form bounded sequences in $L^2(\R^d)$. Set $F_n = E_n U^{(n)}_\lambda f_n$. Then we aim to prove $\|F_n\|_{H^{\alpha/2} (\R^d)} \leq c$ with $c>0$ independent of $n$.   Define $V_n = \{z\in \R^d: |z|_\infty < 2n^{-1} \}$. Moreover, let $z_1^{(n)},\ldots, z_{2^d}^{(n)}$ be the corners of $\mathfrak{Q}_n(0)$.
  We write
  \begin{align}
  \begin{split}
     & \iint\limits_{\R^d\times\R^d} \frac{\big(F_n (\xi)- F_n(\zeta)\big)^2}{ \abs{\xi-\zeta}^{d+\alpha}} d\xi\,d\zeta = \sum_{x\in \scaleZ{n}}\,\,\iint\limits_{\mathfrak{Q}_n(x)\times V_n} \frac{\big(F_n(\xi+\eta)-F_n(\xi)\big)^2}{\abs{\eta}^{d+\alpha}} d\eta\, d\xi \\
        &\quad\quad + \sum_{x\in \scaleZ{n}}\,\,\iint\limits_{\mathfrak{Q}_n(x)\times (\R^d\setminus V_n)} \big(F_n(\xi+\eta)-F_n(\xi)\big)^2 \abs{\eta}^{-d-\alpha} d\eta\, d\xi =: (I_1) + (I_2)
  \end{split}\label{eq:estimateHalphanorm1}
  \end{align}
  Let us first look at $(I_1)$. For $x \in \scaleZ{n}$, $\xi \in \mathfrak{Q}_n(x)$, $\eta \in V_n$, by Taylor's formula
  \begin{align*}
  \frac{\big(F_n(\xi+\eta)-F_n(\xi)\big)^2}{\abs{\eta}^{d+\alpha}} \leq c_0 n^2\sum_{\substack{\widetilde{x}\in \overline{\mathfrak{Q}_n(x)}\cap \scaleZ{n}\\ \widetilde{y}\in \overline{\mathfrak{Q}_n(x)}+\overline{V}_n\cap \scaleZ{n}}} \frac{\big(F_n(\widetilde{x})-F_n(\widetilde{y})\big)^2}{\abs{\eta}^{d+\alpha-2}}
  \end{align*}

  Furthermore,
  \begin{align*}
    \iint\limits_{\mathfrak{Q}_n(x)\times V_n} \abs{\eta}^{2-d-\alpha} d\eta\, d\xi \leq n^{-d} \int\limits_{\abs{\eta}\leq 2\sqrt{d}n^{-1}} \abs{\eta}^{2-d-\alpha} d\eta \leq c_1 n^{-2+\alpha-d}
  \end{align*}

  where $c_1>0$ only depends on $\alpha$ and $d$. Since $\bigabs{\widetilde{y}-\widetilde{x}}\leq 4\sqrt{d}n^{-1}$ the first sum in \eqref{eq:estimateHalphanorm1} can be bounded from above by
  \begin{align*}
    c_2 n^{-2d} &\sum_{x\in \scaleZ{n}} \sum_{\substack{\widetilde{x}\in \overline{\mathfrak{Q}_n(x)}\cap \scaleZ{n}\\ \widetilde{y}\in \overline{\mathfrak{Q}_n(x)}+\overline{V}_n\cap \scaleZ{n}\\\widetilde{x}\neq\widetilde{y}}}\big(F_n(\widetilde{x})-F_n(\widetilde{y})\big)^2 \abs{\widetilde{x}-\widetilde{y}}^{-d-\alpha}\\
    &\leq c_3 n^{-2d}\sum_{\substack{x,h\in \scaleZ{n}\\ 0<|h|_\infty \leq 6n^{-1}}} \big(F_n(x+h)- F_n(x)\big)^2\bigabs{h}^{-d-\alpha}\,.
  \end{align*}
As expected, $(I_1)$ tends to zero for large $n$. In order to tackle $(I_2)$ note that, for all $h,x\in\scaleZ{n}$, $\eta\in \mathfrak{Q}_n(x+h)$, $\xi\in \mathfrak{Q}_n(x)$
  \begin{align*}
    \big(F_n(\eta) -F_n(\xi)\big)^2 &\leq \max_{i,j= 1,\ldots, 2^d} \big(F_n(x+h+z^{(n)}_i)- F_n(x+z^{(n)}_j)\big)^2\\
        &\leq \sum_{i,j}^{2^d} \big(F_n(x+h+z^{(n)}_i)- F_n(x+z^{(n)}_j)\big)^2 \,.
  \end{align*}
For any $\mathfrak{Q}_n(h) \subset \R^d\setminus V_n$, $h\in \Z^d$, and any $\eta \in \mathfrak{Q}_n(h)$ and $i,j= 1,\ldots, 2^d$
  \begin{align*}
    \abs{\eta} \geq |\eta|_\infty \geq \frac12 \bigabs{h + z_i^{(n)} - z_j^{(n)}}_\infty \geq \frac{1}{2\sqrt{d}} \bigabs{h + z_i^{(n)} - z_j^{(n)}} \geq c_4 \bigabs{h + z^{(n)}_i-z^{(n)}_j}\,,
  \end{align*}
where $c_4$ depends only on the dimension. Keeping in mind that the volume of $\mathfrak{Q}_n(x)\times \mathfrak{Q}_n(x+h)$ is $n^{-2d}$ we estimate $(I_2)$ in \eqref{eq:estimateHalphanorm1} from above by
  \begin{align*}
    c_4^{\alpha+d}n^{-2d}&\sum_{i,j=1}^{2^d} \sum_{\substack{x,h\in \scaleZ{n}\\ h + z^{(n)}_i-z^{(n)}_j\neq 0}} \big(F_n(x+h+z^{(n)}_i)- F_n(x+z^{(n)}_j)\big)^2\bigabs{h + z^{(n)}_i-z^{(n)}_j}^{-d-\alpha}\\
    &\leq c_5 n^{-2d}\sum_{\substack{x,h\in \scaleZ{n}\\ h\neq 0}} \big(F_n(x+h)- F_n(x)\big)^2\bigabs{h}^{-d-\alpha}\,.
  \end{align*}

  Hence we obtain by (A3) and \eqref{eq:resolventapprox}
  \begin{align*}
    \iint\limits_{\R^d\times\R^d} &\big(E_n U_\lambda^{(n)}f_n(\xi)- E_n U_\lambda^{(n)}f_n(\zeta)\big)^2 \abs{\xi-\zeta}^{-d-\alpha} d\xi\,d\zeta\\
        &\leq c_5 n^{-2d}\sum_{\substack{x,h\in \scaleZ{n}\\h\neq 0}} \big(U_\lambda^{(n)}f_n(x+h)- U_\lambda^{(n)}f_n(x)\big)^2 \abs{h}^{-d-\alpha} \displaybreak[0]\\
        &\leq c_6 \mathcal{E}^{(n)} \big(U_\lambda^{(n)} f_n, U_\lambda^{(n)} f_n\big)= c_6 \big(U_\lambda^{(n)} f_n, f_n\big) -c_6\lambda\bignorm{U_\lambda^{(n)} f_n}^2\,.
  \end{align*}

  The right-hand side is bounded in $n$. We conclude that $E_n U^{(n)}_\lambda f_n$ is a bounded sequence in the Sobolev space $H^{\alpha/2}(\R^d)$.

  In our situation this means that there is a subsequence $(n''')$ of $(n'')$ such that $E_{n'''} U^{(n''')}_\lambda f_{n'''}$ converges weakly in $H^{\alpha/2}(\R^d)$ and strongly in $L^2(K)$ for $K \subset \R^d$ compact to an element $\widetilde{F}\in H^{\alpha/2}(\R^d)$ for $n''' \to \infty$. Since $E_{n''} U^{(n'')}_\lambda f_{n''}\to U_\lambda f$ pointwise, $U_\lambda f = \widetilde{F}$ almost everywhere. In particular $U_\lambda f\in H^{\alpha/2}(\R^d)$.

  {\bf Step 2:}  Setting for $r\in(0,1)$ $\mathcal{E}_r (f,g) = \frac 12\iint\limits_{\abs{x-y}\geq r} (f(y)-f(x))(g(y)-g(x))k(x,y)\,dx\,dy $ one observes
  \begin{align}\label{eq:ConvergenceStep2}
    \bigabs{\mathcal{E}_r (U_\lambda f, g) - \mathcal{E}(U_\lambda f, g)}\to 0 \quad\text{ for }r\to0.
  \end{align}

  {\bf Step 3:} In analogy to Step 2 set
  \begin{align*}
    \mathcal{E}^{(n)}_r (f_n, g_n) &= \frac{n^{-2d}}{2} \sum_{\substack{x,y\in \scaleZ{n}\\\abs{x-y}\geq r}}(f_n (y)- f_n(x))(g_n(y)-g_n(x)) C^n(x, y)\,.
  \end{align*}

  Then for any $r>0$
  \begin{align}\label{eq:ConvergenceStep3}
    \bigabs{\mathcal{E}^{(n'')}_r (U^{(n'')}_\lambda f_{n''}, g_{n''}) - \mathcal{E}_r (U_\lambda f, g) } \to 0 \quad\text{ for }n'' \to \infty\,.
  \end{align}
Assertion (\ref{eq:ConvergenceStep3}) is easily established by estimating the difference for large enough $n''$ from above by
  \begin{align*}
    \iint\limits_{\abs{x-y}\geq r/2}\Bigabs{\big(U^{(n'')}_\lambda f([y]_{n''})- U^{(n'')}_\lambda f([x]_{n''})\big)\big(g([y]_{n''})-g([x]_{n''})\big)C^{n''}([x]_{n''},[y]_{n''})\\
    \quad\quad - \big( U_\lambda f(y)-  U_\lambda f(x)\big)\big(g(y)-g(x)\big)k(x,y)}\,dx\,dy
  \end{align*}
This term tends to zero for $n'' \to 0$ since $g([x]_{n''}) \to g(x)$ and $U^{(n'')}_\lambda f([x]_{n''})\to U_\lambda f(x)$ uniformly on compacts for $n''\to\infty$. Using (A5) we see that the integrand converges uniformly to $0$. \eqref{eq:ConvergenceStep3} follows by dominated convergence.

  {\bf Step 4:}
  \begin{align}\label{eq:ConvergenceStep4}
    \sup_{n''}\bigabs{\mathcal{E}^{(n'')} \big(U^{(n'')}_\lambda f_{n''}, g_{n^\prime}\big)- \mathcal{E}^{(n'')}_r \big(U^{(n'')}_\lambda f_{n''}, g_{n''}\big)} \to 0\quad\text{ for }r\to 0.
  \end{align}
  Recall that both, $f$ and $g$ have compact support, i.e. the number of elements in $\supp(g) \cap \scaleZ{n}$ is of order $n^2d$. Using Cauchy-Schwartz and (A2) we obtain
\begin{align*}
    \Bigabs{\sum_{\substack{\abs{x-y}<r\\x,y\in \scaleZ{n}}} &\big(U^{(n)}_\lambda f_n(x) - U^{(n)}_\lambda f_n(y)\big)\big(g_n(x)-g_n(y)\big) C^n(x,y)n^{-2d}}^2\\
    &\quad\leq \Big(\sum_{\substack{\abs{x-y}<r\\x,y\in \scaleZ{n}}}\big(U^{(n)}_\lambda f_n(x) - U^{(n)}_\lambda f_n(y)\big)^2C^n(x,y)n^{-2d}\Big) \\
    &\quad\quad\quad\quad \times \sum_{\substack{\abs{x-y}<r\\x,y\in \scaleZ{n}}}\big(g_n(x)-g_n(y)\big)^2 C^n(x,y)n^{-2d}\displaybreak[0]\\
    &\quad\leq c \bigabs{\mathcal{E}^{(n)} \big(U^{(n)}_\lambda f_n, U^{(n)}_\lambda f_n\big)} \, \bigabs{\mathcal{E}_r \big(g_n, g_n\big)} \,.
  \end{align*}
Let us look at the above estimate for $n=n''$. Since $\bigabs{\mathcal{E}^{(n'')} \big(U^{(n'')}_\lambda f_{n''}, U^{(n'')}_\lambda f_{n''}\big)}$ and $\bigabs{\mathcal{E} \big(g_{n''}, g_{n''}\big)}$ are bounded uniformly in $n''$
\[ \supl_{n''} \bigabs{\mathcal{E}^{(n'')} \big(U^{(n'')}_\lambda f_{n''}, U^{(n'')}_\lambda f_{n''}\big)} \, \bigabs{\mathcal{E}_r \big(g_{n''}, g_{n''}\big)} \to 0 \text{ for } r \to 0  \,. \] Step 4 is completed.


  Finally, combining \eqref{eq:ConvergenceStep2},  \eqref{eq:ConvergenceStep3} and \eqref{eq:ConvergenceStep4} by a standard chaining argument proves \eqref{eq:resolventeqlimit}.
\end{proofo}

\section{Approximation of jump processes by Markov chains}\label{sect:approx}

Here, we prove Theorem \ref{theo:approx}.

\begin{proofo}{Proof of Theorem \ref{theo:approx}}
  Obviously, by the symmetry of $k$ (A1) holds. The upper bounds on $k$ imply for $x,y \in \scaleZ{n}$, $\abs{x-y}_\infty \geq 2 n^{-1}$:
  \begin{align*}
    C^n(x,y) &\leq \kappa_5n^{2d}\int\limits_{\substack{\abs{x-\xi}_\infty< n^{-1}/2\\\abs{y-\zeta}_\infty<n^{-1}/2}} \abs{\xi-\zeta}^{-d-\alpha} d\xi\,d\zeta  \\
    &\leq \kappa_5(4d)^{(d+\alpha)/2}n^{2d}\int\limits_{\substack{\abs{x-\xi}_\infty< n^{-1}/2\\\abs{y-\zeta}_\infty<n^{-1}/2}} \abs{x-y}^{-d-\alpha} d\xi\,d\zeta\leq \kappa_5(4d)^{(d+\alpha)/2}\abs{x-y}^{-d-\alpha} .
  \end{align*}

  since $\abs{\xi-\zeta}\geq\abs{\xi-\zeta}_\infty \geq \abs{x-y}_\infty-n^{-1}\geq \frac 12\abs{x-y}_\infty\geq \frac 12 d^{-1/2}\abs{x-y}$. In the same way one shows $C^n(x,y)\geq \kappa_4 (4d)^{-(d+ \alpha)/2}\abs{x-y}^{-d-\alpha}$ for $\abs{x-y}_\infty \geq 2 n^{-1}$. Therefore, (A2), (A3) and (A4) are satisfied.

  Let $k_n (x,y) = C^n([x]_n,[y]_n)$. Fix a compact rectangular subset $\Omega = \prod_{i=1}^{2d} [a_i, b_i]$, $a_i<b_i$ of $\R^d\times\R^d\setminus \diag$ and define $\Omega_\varepsilon = \prod_{i=1}^{2d} [a_i-\varepsilon, b_i+\varepsilon]$ for $\varepsilon>0$. In particular, there is $c_1>0$ with $\leb(\Omega_\varepsilon\setminus\Omega)<c_1\varepsilon$ for $\varepsilon<1$. Then there exists $n_0>0$ such that $\dist (\Omega, \diag)> n_0^{-1}$. For all $n>4n_0$ and $\varepsilon<(4n_0)^{-1}$ the functions $k_n$ are uniformly bounded from above on $\Omega_\varepsilon$ by $c_2>0$. By the Theorem of Lusin, for each $\varepsilon>0$ there exists a compact set $K_\varepsilon\subset\Omega$ such that $k$ restricted to $K_\varepsilon$ is continuous while $\leb(\Omega\setminus K_\varepsilon)<\varepsilon$. Furthermore, there exists a continuous function $k^\varepsilon\colon\R^d\times\R^d\to\R^+$ with compact support such that $k^\varepsilon=k$ on $K_\varepsilon$, $\supp k^\varepsilon\subset\Omega_{\varepsilon}$ and $\onorm{k^\varepsilon}_\infty < c_2$. We estimate
  \begin{align*}
    \onorm{k-k^\varepsilon}_{L^1(\Omega)} &=\int_\Omega \oabs{k-k^\varepsilon} = \int_{K\varepsilon} \oabs{k-k^\varepsilon}+\int_{\Omega\setminus K_\varepsilon} \oabs{k-k^\varepsilon}\leq c_2\varepsilon.
  \end{align*}
   Define now $k_n^\varepsilon(x,y) = C_\varepsilon^n ([x]_n, [y]_n)$ as above with $k$ replaced by $k^\varepsilon$. Then, since $k^\varepsilon$ is Riemann integrable with compact support, $k^\varepsilon_n$ converges in $L^1(\R^{2d})$ for $n\to \infty$ to $k^\varepsilon$. Moreover, for $\varepsilon<(4n_0)^{-1}$, $n>\varepsilon^{-1}$ and by the definition of the conductivity functions
   \begin{align*}
     \onorm{k_n^\varepsilon-k_n}_{L^1(\Omega)} &\leq \onorm{k^\varepsilon-k}_{L^1(\Omega_\varepsilon)} \leq \onorm{k^\varepsilon-k}_{L^1(\Omega)} + \onorm{ k^\varepsilon-k}_{L^1(\Omega_\varepsilon-\Omega)} \leq c_2(1+c_1) \varepsilon
   \end{align*}
   Putting all this together we get for $n$ large enough
   \begin{align*}
     \onorm{k-k_n}_{L^1(\Omega)} &\leq \onorm{k-k^\varepsilon}_{L^1(\Omega)} +\onorm{k^\varepsilon-k_n^\varepsilon}_{L^1(\Omega)} +\onorm{k_n^\varepsilon-k_n}_{L^1(\Omega)}\leq c_3\varepsilon.
   \end{align*}
   This directly yields (A5). Now Theorem \ref{theo:meta-one} applies.
\end{proofo}

\def\cprime{$'$}


{\bf Ryad Husseini}\\
Institut f\"{u}r Angewandte Mathematik\\
Universit\"{a}t Bonn\\
Poppelsdorfer Allee 82\\
D-53115 Bonn, Germany\\
ryad@uni-bonn.de\\

{\bf Moritz Kassmann}\\
Institut f\"{u}r Angewandte Mathematik\\
Universit\"{a}t Bonn\\
Beringstrasse 6\\
D-53115 Bonn, Germany\\
kassmann@iam.uni-bonn.de

\end{document}